\pdfoutput=1
\documentclass[10pt]{article}

\usepackage[left=1.25cm,right=1.25cm,top=2cm,bottom=2cm]{geometry}
\usepackage[width=0.9\textwidth]{caption}

\usepackage[titletoc, title]{appendix}

\usepackage{amsfonts,amssymb,amstext,amsmath,cases}
\usepackage{amsthm}
\usepackage[overload]{empheq}%
\usepackage{multirow}

\usepackage{mathtools}

\usepackage[none]{hyphenat}
\hyphenchar\font=-1
\usepackage{enumitem}
\usepackage[colorlinks,citecolor=blue,linkcolor=blue,urlcolor=blue]{hyperref}
\usepackage[hyperpageref]{backref}

\usepackage{multirow}
\usepackage{array}
\usepackage{booktabs}
\usepackage{makecell}

\usepackage{comment}

\usepackage{graphicx}
\usepackage{caption}

\usepackage{xcolor}
\usepackage{subcaption}

\newtheorem{remark}{Remark}

\newcommand{\bd}[1]{\boldsymbol{#1}}
\newcommand{\wtd}[1]{\widetilde{#1}}
\newcommand{\wh}[1]{\widehat{#1}}
\newcommand{\abs}[1]{\left\vert #1 \right\vert}
\newcommand{\ds}{\displaystyle}
\newcommand{\dint}{\displaystyle\int}
\newcommand{\dsum}{\displaystyle\sum}
\renewcommand{\r}{\rho}
\newcommand{\vv}{\bd{v}}
\newcommand{\p}{\partial}

\graphicspath{{figuras/}{./}}

\title{General function bases moment models for suspended sediment with vertical profile of concentration}

\author{C. Caballero-Cárdenas\thanks{Facultad de Negocios y Tecnología, Universidad Alfonso X el Sabio - Mare Nostrum. Camino de la Térmica, 90. 29004 Málaga, Spain ({\tt ccabacar@uax.es}).} \and J. Garres-Díaz \thanks{$^{\mathrm{C}}$Departamento de Matemática Aplicada II $\&$ IMUS, ETS Ingeniería, Universidad de Sevilla, Camino de los Descubrimientos s/n, 41092 Sevilla, Spain {\tt (jgarres@us.es)}}\, \and
I. Gómez-Bueno\thanks{Dpto. Matem{\'a}tica Aplicada. Universidad de M\'alaga, 29071 M\'alaga, Spain 
	({\tt igomezbueno@uma.es}).} 
\and T. Morales de Luna\thanks{Dpto. An\'alisis Matem\'atico, Estad. e I.O. y Matem\'atica Aplicada. Universidad de M\'alaga, 29071 M\'alaga, Spain  ({\tt tmorales@uma.es})}
}

\begin{document}
\date{}
\maketitle

\begin{abstract}
	We introduce a general hydrostatic shallow moment framework for suspended sediment transport with vertically varying horizontal velocity and concentration profiles. The formulation is constructed over arbitrary vertical functional bases that are numerically orthogonalised via the Gram-Schmidt algorithm. This approach allows the approximation spaces for the vertical profile to be adapted to the specific structure of each physical variable, where also the velocity and sediment concentration can be projected onto different functional spaces. To overcome the limitations due to the inherent oscillatory nature of global Legendre polynomials, we explore a specific continuous piecewise linear $\mathbb{P}_1$ finite-element basis. Furthermore, non-uniform vertical partitions are considered, enabling an enhanced representation of the vertical profiles in specific regions. Numerical experiments involving boundary layers and suspended sediment transport demonstrate that the proposed framework is accurate to reproduce the expected vertical profiles. The improvement is particularly significant for sediment concentration profiles with sharp vertical variations, where polynomial reconstructions exhibit pronounced oscillations leading to physically meaningless profiles, whereas the $\mathbb{P}_1$ bases produce a qualitatively accurate representation of the vertical profiles observed in the sedimentation processes.\\
	
\end{abstract}

\section{Introduction}\label{se:intro}

The mathematical modeling and numerical simulation of geophysical fluid dynamics play a crucial role in understanding environmental processes, predicting natural hazards, and managing coastal and river engineering systems. Among these processes, sediment transport and its interaction with the underlying fluid flow are of primary importance due to their profound impact on morphodynamics, water quality, and infrastructure safety. Sediment transport is generally categorized into bedload transport, where grains roll or slide along the movable bottom modifying the topography, and suspended sediment transport, where fine particles are lifted and carried by currents within the fluid column (see, e.g., \cite{Dey2014}). In this work we focus on suspended sediment transport, which introduces highly stratified flows where concentration varies significantly across the vertical water column. 

Accurately simulating suspended sediment requires models capable of reproducing these sharp vertical distributions without incurring huge computational costs. A widely known strategy to simulate free-surface flows efficiently is the classical shallow water system (see e.g. \cite{SaintVenant1871, Gerbeau2001}). It is based on integrating the Navier-Stokes equations over the fluid depth under the hydrostatic assumption, so that the system's unknowns are the water depth and the depth-averaged horizontal velocity. However, a major limitation of standard shallow water models is the loss of vertical variations, since they assume a vertically uniform profile of the horizontal velocity. This assumption fails dramatically in many geophysical applications. Concretely, in our specific case of study, sediment particles tend to concentrate near the bottom, exhibiting steep vertical profiles that are far from uniform, leading to a significant deviation of the horizontal velocity from a constant profile. Nevertheless, in many classical studies, the concentration profile is usually assumed to be vertically uniform for simplicity (see e.g. \cite{parker1986self,GonzalezAguirre2020,bradford1999hydrodynamics,khan2005numerical, Burger2025}).

To overcome this limitation of depth-averaged models, several strategies to recover vertical effects have been proposed in the literature, notably multilayer models \cite{FernandezNieto2013a} (including applications to suspended sediment \cite{FernandezNieto2013, MoralesdeLuna2017, Buerger2019}) and moment-based approaches \cite{Steffler1993, Kowalski2018} (with application to suspended sediment transport as in \cite{Guo1999,CanteroChinchilla2019}). A comprehensive comparative study between these two frameworks was recently addressed in \cite{GarresDiaz2023}, where a family of \textit{multilayer-moment} models combining both strategies was introduced. The core idea of the multilayer approach is to consider a vertical partition of the fluid domain and assume the shallowness hypothesis within each layer, thereby obtaining a piecewise constant reconstruction of the variables. On the contrary, the moment approach expands the vertical dependency of the variables within a functional space, where the corresponding expansion coefficients (\textit{moments}) become the unknowns of the system. In most works dealing with shallow moment models, the orthogonal scaled Legendre polynomials are chosen as functional basis (see e.g. \cite{Kowalski2018, CanteroChinchilla2019, koellermeier2020analysis, Koellermeier2022}). Recent developments have also focused on the mathematical structure of these models, including their rotational invariance and hyperbolicity in two spatial dimensions \cite{Bauerle2025}. One of the main conclusions in \cite{GarresDiaz2023} is that moment models are computationally more efficient than multilayer methods, although the latter offer greater flexibility to reproduce arbitrary vertical profiles. Conversely, moment models based on polynomials fail to accurately reproduce profiles containing almost constant regions, such as velocity in fluids with boundary layer at the bottom. This drawback is due to the inherently oscillatory nature of global polynomials. Interestingly, a recent work \cite{Koellermeier2025} considered linear and quadratic splines as an alternative vertical basis. Nevertheless, the resulting basis was non-orthogonal and then a full coefficient mass matrix had to be inverted, which increased the complexity of the final system. To our knowledge, shallow moment models derived over an arbitrary functional basis have not yet been addressed in the literature, although such an approach would provide a flexible and general framework easily adaptable to the specific vertical profiles of any given application.

Concerning suspended sediment transport, on the one hand, several multilayer models for polydisperse sedimentation capable of correctly reproducing vertical concentration profiles have been introduced \cite{FernandezNieto2013, MoralesdeLuna2017, Buerger2019}. On the other hand, although shallow moment systems have already been developed for bedload transport problems (see \cite{CanteroChinchilla2019, GarresDiaz2021b}), and \cite{Guo1999} considered a model with a linear expansion of both the horizontal velocity and sediment concentration, the derivation of a moment-based shallow water model considering the vertical stratification of suspended sediment remains unaddressed in the literature.

In this paper, we introduce a general derivation of a hydrostatic moment model for suspended sediment transport where both horizontal velocity and sediment concentration are expanded in an arbitrary functional basis. As a first step, once the basis is chosen, the Gram-Schmidt orthogonalisation algorithm is numerically applied to construct an orthogonal basis, which notably simplifies the computations, the final form of the system, and consequently the numerical methods used to approximate it. Thus, the contribution of this framework is twofold: (i) to our knowledge, this is the first shallow moment model dealing with a vertical expansion of the sediment concentration; and (ii) even for clear water flows, for which our proposal remains fully valid as a standard shallow water moment model (with no sediments), the derivation of the model for an arbitrary basis already constitutes a significant contribution. Furthermore, the resulting model is highly flexible since each variable (velocity and concentration) may be expanded in a different functional subspace. 

Concretely, to avoid the oscillatory nature of Legendre polynomials, we explore a specific basis consisting of continuous piecewise linear polynomials $\mathbb{P}_1$. In addition, we consider non-uniform vertical partitions, thereby achieving an enhanced vertical resolution in specific regions of interest. In the numerical experiments, this approach is applied to flows involving boundary layers and suspended sediment. The vertical profiles of horizontal velocity and sediment concentration obtained notably improve the results of the polynomial case. In particular, a qualitatively accurate representation of the sedimentation process is achieved.

The paper is organized as follows. In Section \ref{se:derivationModel}, we present the derivation of the hydrostatic shallow moment model for suspended sediments over an arbitrary orthogonal basis. Section \ref{se:finalModel} summarizes the final model, including an asymmetric formulation where the velocity and sediment concentration are projected onto different functional spaces. In addition, Subsection \ref{se:ParticularBasis} introduces the piecewise linear bases over a non-uniform partition, which is subsequently used in the numerical experiments. Section \ref{se:numericalApprox} is devoted to the numerical approximation of the final system within the finite volume framework. Several numerical tests concerning water flows and sedimentation processes are presented in Section \ref{se:numericalTest}. In particular, the resulting profiles obtained with the proposed functional bases significantly improve upon the results of the Legendre polynomial bases. Finally, some concluding remarks are in Section \ref{se:conclusions}.

\section{Derivation of a hydrostatic moment model for a general functional basis}\label{se:derivationModel}
In this section, we derive the hydrostatic moment model for suspended sediments with a vertical concentration profile. We begin by describing the physical system of equations governing the evolution of the sediment concentration. Next, the moment system is derived considering a general orthogonal functional basis. Finally, we explore a particular choice of functional basis that allows us to describe the vertical profile of sediment concentrations, in contrast to the conventional polynomial bases used in the literature. For the sake of simplicity, we consider here a single species of sediment and the 1D case, although the approach could be generalized. 

\subsection{Governing equations}\label{se:sistemaInicial}
Let us consider a sediment species with density $\r_s$ and velocity $\vv_s=(u_s,w_s)^T$ inside a fluid with density $\r_f$ and velocity $\vv_f=(u_f,w_f)^T$. Denoting the sediment concentration by $c\in[0,1)$, the continuity equations for each phase read:
\begin{subequations}
	\label{eq:continuity}
	\begin{eqnarray}
		&&    \r_s\left(\p_t c + \nabla\cdot(c \vv_s)\right) = 0,\label{eq:continuity_a}\\
		&&    \r_f\left(\p_t (1-c) + \nabla\cdot((1-c) \vv_f)\right) = 0.\label{eq:continuity_b}
	\end{eqnarray}
\end{subequations}
From the previous equations, we obtain:
\begin{equation}
	\label{eq:divV}
	\nabla\cdot \vv = 0 \qquad\mbox{with}\qquad \vv = c\vv_s + (1-c)\vv_f,
\end{equation}
where $\vv=(u,w)^T$ is the volume-averaged mixture velocity. Now, we assume that there is no horizontal slip between the phases, while vertical slip may occur. Then, each phase velocity ($\vv_s$ and $\vv_f$) can be written in terms of its deviation $(\delta \vv)$ from the average velocity $\vv$ as:
\begin{equation*}
	\vv_k = \vv + (\delta \vv)_k \quad\mbox{with}\quad (\delta \vv)_k = (0,(\delta w)_k)^T, \qquad \qquad \mbox{for }k=f,s.
\end{equation*}
Using the previous decomposition, the continuity equation \eqref{eq:continuity_a} can be rewritten in terms of the average velocity as:
$$
\p_t c +  \nabla\cdot(c \vv) + \p_z \left(c (\delta w)_s\right) = 0,
$$
where the term $(\delta w)_s$ is the vertical slip velocity, which will be defined through an empirical law later. 

Notice that the average velocity may also be written as $\vv = \vv_f  + c(\vv_s-\vv_f) = \vv + (\delta\vv)_f + c((\delta\vv)_s-(\delta\vv)_f)$, which allows us to express:
$$
(\delta \vv)_f = \dfrac{-c}{1-c} (\delta \vv)_s.
$$

Additionally, as a hypothesis, we adopt the Boussinesq approximation, which states that in the case of variable density, the deviations from a constant reference value are small and their effects only appear in the pressure and gravitational terms. To this end, we write the mixture density $\r$ as:
$$
\dfrac{\r}{\r_f} = \dfrac{\r_s c + \r_f(1-c)}{\r_f} = 1 + c\Delta \r, \quad\mbox{with}\quad \Delta\r = \dfrac{\r_s-\r_f}{\r_f}, 
$$
where $c\Delta \r  \sim \mathcal{O}(10^{-1})$ for low concentrations. Therefore, as a starting system, we consider the incompressible Euler equations for the fluid-sediment mixture, together with the previous equation for the concentration ($c$) and the Boussinesq approximation, namely:
\begin{subequations}\label{eq:initialSystem}
\begin{empheq}[left={\empheqlbrace\,}]{align}
	&\p_t c + \nabla\cdot(c \vv)+ \p_z \left(c (\delta w)_s\right) = 0,\label{eq:initialSystem_a}
	\\[2mm]
	&\nabla\cdot \vv=0,\label{eq:initialSystem_b}
	\\
	&\partial_t \vv+\nabla\cdot(\vv\otimes \vv)+\frac{1}{\rho_f}\nabla p
	=
	\frac{\rho}{\rho_f}\bd{g},\label{eq:initialSystem_c}
\end{empheq}
\end{subequations}
where $\bd{g}=(0,-g)^T$ is the gravitational acceleration, and $p$ is the pressure. In this work we shall assume the pressure to be hydrostatic. Therefore, convective effects are neglected in the vertical momentum equation, which means that we set:
\begin{equation}
\label{eq:presHidro}
\dfrac{1}{\r_f}\p_z p = -(1+c\Delta\r)g.
\end{equation}
As boundary conditions we consider the usual kinematic conditions at the bottom and free surface and no-stress at the free surface, that is,
\begin{equation}
	\label{eq:condContorno}
	\p_t\eta + \vv_{|_{z=\eta}}\cdot \left(\p_x\eta,-1\right) = 0, \qquad  \vv_{|_{z=b}}\cdot \left(\p_x b,-1\right) = 0, \qquad p_{|_{z=\eta}} = 0,
\end{equation}
where $\eta$ and $b$ denote the free surface and bottom levels.
\begin{remark}
The mixture velocity $\vv$ (see \eqref{eq:divV}) used above is not the barycentric velocity of the mixture, which is defined as
\begin{equation}
	\widehat{\vv}\equiv(\widehat{u},\widehat{w})^T=\dfrac{1}{\rho}\left(\r_s c \vv_s+\r_f (1-c) \vv_f\right),\qquad \mbox{with}\qquad \r = \r_s c +\r_f (1-c).
\end{equation}    
However, one can write:
$$
\widehat{\vv} = \vv + c \Delta\r\dfrac{\rho_f}{\rho}(\delta \vv)_s,\quad \mbox{and then}\quad \left\{\begin{array}{l}
	\widehat{u} = u,\\
	\widehat{w} = w + c \Delta\r\dfrac{\rho_f}{\rho}(\delta w)_s.
\end{array}\right.
$$
where $\r_f/\r \sim \mathcal{O}(1)$, the slip velocity is typically small, and $c\Delta \r  \sim \mathcal{O}(10^{-1})$ under the Boussinesq approximation. Consequently, the momentum equations \eqref{eq:initialSystem_c} written in terms of $\widehat{\vv}$ and $\vv$ coincide when the Boussinesq approximation is applied; that is, the deviation from a constant density only affects the pressure terms. 
\end{remark}

\subsubsection*{Closure relation for the vertical slip}
A closure relation defining the vertical slip velocity $(\delta w)_s$ is required. Following \cite{MoralesdeLuna2017}, the following empirical law, adapted to the case of a single sediment species, is adopted:
\begin{subequations}
\label{eq:deltaW}
\begin{equation}
	(\delta w)_s = -\omega_{s} \, \chi(c),
\end{equation}
with $\omega_s$ a constant value known as the settling velocity, which depends on the type of sediments considered, and the function $\chi(c)$, which models the interaction effects between particles, defined as:
\begin{equation}
	\chi(c) = (1-c)^m\,\mathbf 1_{\{c<c_{\max}\}} = 
	\begin{cases}
		\left(1 - c \right)^m,
		& \text{if } c < c_{\max}, \\[0.4em]
		0, & \text{otherwise}.
	\end{cases}
\end{equation}
\end{subequations}
Here, $m\in[2.5,5.5]$ is a fixed parameter, and $c_{\max}$ is the saturation concentration. This definition of $\chi(c)$ implies that the effective settling velocity decreases as the sediment concentration increases, becoming zero when the maximum threshold $c_{\max}$ is reached. Thus, in practice, the concentration satisfies $0<c(t,x,z)<c_{\max}$, where a typical value for $c_{\max}$ is $0.8$.

\subsection{Moment approach in orthogonal bases}
In this section, we perform the vertical moment discretization of the system \eqref{eq:initialSystem}, \eqref{eq:presHidro}, \eqref{eq:deltaW} describing a sediment-fluid mixture with velocity $\vv = (u,w)^T$ and sediment concentration $c(t,x,z)$.  As usual, the first step consists of writing the initial system in $\sigma$-coordinates. After, the equations of the resulting system are multiplied by basis functions and integrated to get the final shallow moment model. In the following, this procedure is developed. 
  
\subsubsection{$\sigma$-coordinates}
Let us recall the main ingredients of the $\sigma$-coordinates change of variables. Defining the total depth $h(t,x)$ and a fixed bottom topography $b(x)$, the change of variable 
\begin{equation}
	\label{eq:cambioVariable}
	\xi = \frac{z - b}{h}
\end{equation}
transforms the moving vertical domain $z\in[b, b+h]$ into the constant interval $\xi\in[0,1]$. Then, for a function $\psi(t,x,z)$, we have the corresponding definition in $\sigma$-coordinates $\tilde{\psi}(t,x,\xi) := \psi(t,x,\xi h + b)$, together with the differential operators:
$$
\partial_\xi \tilde{\psi} = h\,\partial_z \psi, \qquad \mbox{and}\qquad     h\,\partial_s \psi=\partial_s (h \tilde{\psi})-\partial_\xi \bigl( \partial_s (\xi h + b)\,\tilde{\psi} \bigr) \quad\mbox{for }s \in \{t,x\}.
$$

Introducing the coupling term
\begin{equation}
\omega(\xi)=\partial_t(\xi h+b)+\tilde{u}\,\partial_x(\xi h+b)-\tilde{w},
\end{equation}
the boundary conditions at the free surface ($\xi=1$) and the bottom ($\xi=0$) are 
\begin{equation}
	\label{eq:kinematic}
	\omega(0) = \omega(1) = 0,\qquad \mbox{and}\qquad \wtd{p}(1)=0,
\end{equation} 
and equations \eqref{eq:initialSystem}, \eqref{eq:presHidro}, \eqref{eq:deltaW} take the form:
\begin{subequations}\label{eq:systemSigma}
	\begin{empheq}[left={\empheqlbrace\,}]{align}
		&\partial_t (h \tilde{c}) + \partial_x (h \tilde{c} \tilde{u}) - \partial_\xi \left( \omega \tilde{c} - (\delta w)_s \tilde{c} \right) = 0,\label{eq:systemSigma_a}\\[2mm]
		&\partial_t h + \partial_x (h \tilde{u}) - \partial_\xi \omega = 0, \label{eq:systemSigma_b}
		\\[2mm]
		&\partial_t (h \tilde{u}) + \partial_x \left( h \tilde{u}^2 + \frac{1}{\rho_f}h\tilde{p} \right)
		- \partial_\xi \left( \omega \tilde{u} + \frac{1}{\rho_f}\tilde{p}\,\partial_x(\xi h+b) \right) = 0.\label{eq:systemSigma_c}
\end{empheq}
\end{subequations}
Additionally, the hydrostatic pressure \eqref{eq:presHidro} is written as:
\begin{equation}\label{eq:systemSigma_d}
\dfrac{1}{\r_f}\p_{\xi} \tilde{p} = -(1+\tilde{c}\Delta\r)gh.
\end{equation}
Using the second equation in the previous system, the coupling term can be written as:
\begin{equation}\label{eq:systemSigma_e}
\omega = \partial_x \left( h \left(\int_0^\xi \tilde{u}(s)\,ds - \xi \bar{u}\right)\right), \qquad\mbox{with}\qquad \bar{u} = \ds\int_0^1 \tilde{u}(\xi) d\xi.
\end{equation}

Conventional hydrostatic shallow water moment models consider an expansion of the horizontal velocity
$$\tilde{u}(t,x,\xi) = \sum_{l=0}^{N} u_l(t,x) \phi_l(\xi),$$
within the functional space spanned by the basis $\{\phi_0,\phi_1,\dots,\phi_N\}$, where the coefficient of $\tilde{u}(t,x,\xi)$ in this basis, i.e., $u_l(t,x)$, are the so called \textit{moments}. In many studies, the orthogonal bases given by scaled Legendre polynomials has been used (see, e.g., \cite{Kowalski2018,GarresDiaz2021b,Koellermeier2022,GarresDiaz2023}), although any functional basis could be adopted in principle. For instance, in \cite{Koellermeier2025}, the authors consider linear and quadratic splines as a basis (including $\phi_0=1$) that is not orthogonal. This requires the resulting model to compute the inverse of a coefficient matrix, which, even if computed only once, leads to a more complex system structure, particularly regarding the non-conservative product terms.

A contribution of this paper is the derivation of hydrostatic moment models over a general basis aimed at applications to fluids with suspended sediments, although many other configurations are possible. The core idea is adapting the functional basis $\{\phi_0,\phi_1,\dots,\phi_N\}$ to the expected vertical profiles of velocity for each specific application. Therefore, it is convenient to choose a basis satisfying: (i) $\phi_0 = 1$ to recover $u_0 = \bar{u}$, and (ii) orthogonality to simplify computations and avoid matrix inversions that lead to more complex systems, as mentioned above. Once the basis is chosen, the moment based hydrostatic model is obtained by the usual integration procedure. 

\subsubsection{Integration procedure}\label{sec:integration_proc}

Let $\{\varphi_0(\xi),\varphi_1(\xi),\varphi_2(\xi),\dots,\varphi_N(\xi)\}$ be a set of linearly independent functions in $L^2(0,1)$ which span the corresponding subspace. By applying the Gram-Schmidt orthogonalisation process to the previous set, we obtain an orthogonal basis in the subspace denoted by $\{\phi_0(\xi),\phi_1(\xi),\dots,\phi_N(\xi)\}$, so that 
\begin{equation}
	\int_0^1 \phi_i(\xi)\phi_j(\xi)\,d\xi = \mu_i \delta_{ij}, \qquad \text{where} \quad \mu_j = \int_0^1 \phi_j(\xi)^2\,d\xi.
\end{equation}
In this work we shall assume that the first element of the basis is the constant unitary function, $\varphi_0(\xi)=1,$ and consequently we trivially get $\phi_0(\xi) = 1$. This will allow us to recover the averaged values of the variables directly, as it will be stated afterwards. Now, both velocity and concentration fields are decomposed as
\begin{equation}
	\tilde{u}(t,x,\xi)=\sum_{j=0}^{N}u_j(t,x)\phi_j(\xi),
	\qquad
	\tilde{c}(t,x,\xi)=\sum_{i=0}^{N}c_i(t,x)\phi_i(\xi).
\end{equation}

As a consequence of the orthogonality, $\int_0^1 \phi_k d\xi = 0$ holds for $k>0$. Using that, the mean velocity coincide with the first moment $u_0(t,x)$ since
$$
	\bar{u}(t,x) =  \int_0^1 \tilde{u}(t,x,\xi) d\xi = \sum_{l=0}^N \int_0^1 u_l(t,x)\phi_l(\xi) d\xi  = u_0(t,x),
$$
and, analogously, the mean concentration is $\bar{c}(t,x) =\int_0^1 \tilde{c}(t,x,\xi) d\xi =  c_0(t,x)$. 

\medskip

Before detailing the integration process, let us introduce all the coefficients depending on the basis $\{\phi_0(\xi),\phi_1(\xi),\dots,\phi_N(\xi)\}$ that will appear in the model derivation. We shall remark that, since the function basis is fixed at the beginning, the orthogonalisation process is carried out once and all the required coefficients are saved via a pre-processing step, so that they will be imported and read whenever necessary. So, it does not entail any additional computational cost in the simulations. For  $i,j,k \in \{0, \dots, N\}$, we define
\begin{subequations}
	\label{eq:coeficientes}
\begin{equation}\label{eq:coeficientes1}
	\begin{array}{l}
	A_{ijk} = \dint_0^1 \phi_i\,\phi_j\,\phi_k\, d\xi, \qquad
	B_{ijk} = \dint_0^1 \phi_i'\left(\int_0^\xi \phi_j(s)\,ds\right)\phi_k\, d\xi, \qquad
	C_{ijk} = \dint_0^1 \phi_i'\,\phi_j\,\phi_k\, d\xi,
		\end{array}
\end{equation}
as well as the coefficients
\begin{equation}\label{eq:coeficientes2}
	\begin{array}{l}
		D_{ij} = \dint_0^1 \phi_i\left(\int_\xi^1 \phi_j(s)\,ds\right)\, d\xi, \qquad
		E_{ij} = \dint_0^1 \xi\,\phi_i\,\phi_j\, d\xi,%
	\end{array}
\end{equation}
\end{subequations}
which will appear in the hydrostatic baroclinic pressure depending of the sediment concentration. 

\medskip

The depth-integrated mass conservation equation is immediately obtained by integrating \eqref{eq:systemSigma_b} and using kinematic conditions \eqref{eq:kinematic}:
\begin{equation}\label{eq:ecuacionMasa}
	\partial_t h + \partial_x (h u_0) = 0.
\end{equation}

In order to obtain evolution equations for the moment coefficient of the velocity ($u_k(t,x)$) and concentration ($c_k(t,x)$), equations \eqref{eq:systemSigma_c} and \eqref{eq:systemSigma_a} are multiplied by the function basis $\phi_k$ and integrated along the vertical direction $\xi\in[0,1]$. This is a standard procedure, so we refer the reader to \cite{GarresDiaz2023} for the details and only the computations that are a novelty with respect to previous works will be detailed. Concerning the momentum equation \eqref{eq:systemSigma_c} the main difficulty is related to the pressure terms. 

First, equation \eqref{eq:systemSigma_d} is integrated at the interval $[\xi,1]$, obtaining
\begin{equation}
	\dfrac{1}{\r_f}\tilde{p}(\xi) = g h (1-\xi) + gh \Delta \rho \int_\xi^1 \tilde{c}(\xi) \, d\xi.
\end{equation}
where boundary condition \eqref{eq:kinematic} has been used. The first term on the right represents the barotropic hydrostatic component (due to the weight of clean water), while the second term represents the baroclinic contribution induced by density variations due to the presence of suspended sediments. The latter is responsible for generating turbidity currents and modifying the vertical velocity profile through coupling in the momentum equation.

The pressure terms in \eqref{eq:systemSigma_c} are
$$
\frac{1}{\rho_f}\left[ \partial_x \left(h\tilde{p} \right)
- \partial_\xi \left( \tilde{p}\,\partial_x(\xi h+b) \right)\right] = gh\p_x(b+h) + gh\Delta\r \left[\p_x\left(h\dint_{\xi}^1 \wtd{c}(\tau)d\tau\right) + \wtd{c}(\xi)\left(\xi\p_xh+\p_xb\right)\right].
$$
This expression is multiplied by $\phi_k$ and integrated for $\xi\in[0,1]$ to obtain
\begin{equation}
gh \partial_x (b+h) \delta_{k0} + g \Delta \rho  \left[h\dsum_{i=0}^{N} D_{ki}\p_x(hc_i) + \partial_x h  \sum_{i=0}^N E_{ki}\, (hc_i) + \partial_x b  \sum_{i=0}^N A_{ki0}\, (hc_i)\right].
\end{equation}

Thus, the evolution equations for the velocity modes $u_k(t,x)$ are, for $k=0,\dots,N$,
\begin{multline}\label{eq:finalMomentum_eng}
	 \mu_k\partial_t (h u_k) + \partial_x \left( h \sum_{i,j=0}^N u_i u_j A_{ijk} \right) + \sum_{i,j=0}^N u_i \partial_x (h u_j) B_{kji} - \partial_x (h u_0) \sum_{i=0}^N u_i B_{k0i} \\ 
	+gh \partial_x (b+h) \delta_{k0} + g \Delta \rho  \left[h\dsum_{i=0}^{N} D_{ki}\p_x(hc_i) + \partial_x h  \sum_{i=0}^N E_{ki}\, (hc_i) + \partial_x b  \sum_{i=0}^N A_{ki0}\, (hc_i)\right] = 0.
\end{multline}
Note that the barotropic pressure gradient is filtered out via the Kronecker delta $\delta_{k0}$ for all higher-order moments ($k>0$). Consequently, the free-surface gradient directly affects only the mean velocity ($u_0$), whereas the baroclinic interactions dynamically adapt the vertical velocity distribution based on the suspended sediment concentration profile.

Finally, equation \eqref{eq:systemSigma_a} is considered to obtain the evolution equations for the moments coefficient ($c_k(t,x)$) of the concentration variable. Here, the main difficulty is related to the vertical slip velocity $(\delta w)_s$. Multiplying that equation by $\phi_k$ and integrating we get the evolution equation
\begin{equation}\label{eq:ck1}
		 \mu_k\partial_t (h c_{k}) + \partial_x \left(h \sum_{i,j=0}^N c_i u_j A_{ijk}\right) + \sum_{i,j=0}^N c_i \partial_x (h u_j) B_{kji} - \partial_x(h u_0) \sum_{i=0}^N c_i B_{k0i} + \mathcal S_k(c)
		=0.
\end{equation}   
where $\mathcal S_k(c)$ denotes the contribution of the relative vertical settling velocity. This term deserves special attention. According to \eqref{eq:deltaW}, we define the vertical flux
\begin{equation}
	F(c) = -\omega_{s} c\,\chi(c),
\end{equation}
and then $S_k(c)$ is the projection of the derivative of this vertical flux
\begin{equation}
	\label{eq:sedimentacionMomentos}
\mathcal S_k(\mathbf c) =
\int_0^1 \phi_k(\xi) \partial_\xi F(\tilde c(\xi))\,d\xi .
\end{equation}
In the numerical implementation, this term is not explicitly computed in the shallow moment framework but it is approximated through a reconstruction--evolution--projection procedure on a finite difference vertical grid. It will be detailed later in Section \ref{se:numericalApprox}. 

\section{Final hydrostatic shallow moment model}\label{se:finalModel}
The final hydrostatic shallow moment system for a single sediment species, for an arbitrary orthogonal functional basis $\{\phi_0,\phi_1,\dots,\phi_M\}$ with $\phi_0=1$, for $k=0,...,N$, reads:
\begin{subequations}\label{eq:finalModel}
	\begin{empheq}[left={\empheqlbrace\,}]{align}
		&\partial_t h + \partial_x (hu_0) = 0,\label{eq:finalModel_c1}\\[2mm]
		& \mu_k\partial_t (h u_k) + \partial_x \left( h \sum_{i,j=0}^N u_i u_j A_{ijk} \right) + \sum_{i,j=0}^N u_i \partial_x (h u_j) B_{kji} - \partial_x (h u_0) \sum_{i=0}^N u_i B_{k0i} \nonumber\\ 
		&\qquad+gh \partial_x (b+h) \delta_{k0} + g \Delta \rho  \left[h\dsum_{i=0}^{N} D_{ki}\p_x(hc_i) + \partial_x h  \sum_{i=0}^N E_{ki}\, (hc_i) + \partial_x b  \sum_{i=0}^N A_{ki0}\, (hc_i)\right] = 0, \label{eq:finalModel_u}
		\\[2mm]
		& \mu_k\partial_t (h c_{k}) + \partial_x \left(h \sum_{i,j=0}^N c_i u_j A_{ijk}\right) + \sum_{i,j=0}^N c_i \partial_x (h u_j) B_{kji} - \partial_x(h u_0) \sum_{i=0}^N c_i B_{k0i} + \mathcal S_k(c)
		=0.\label{eq:finalModel_c}
	\end{empheq}
\end{subequations}
Note that in the last term in \eqref{eq:finalModel_u} one can use that $\dsum_{i=0}^N A_{ki0} (hc_i) = hc_k \mu_k$ since $A_{ki0} = \delta_{ki}\mu_k$.

\subsection{Asymmetric Formulation ($N \neq M$)}\label{se:asimetrico}

In applications where the vertical profile of sediment concentration presents sharper gradients than the velocity profile (for example, near the bed), it will be efficient to decouple the order of expansions. 

In principle, one could consider different functional bases for the velocity and concentration expansions. However, the expression of the coefficients in the system would be more complex. For the sake of completeness, we write here the general form of the system, although it will not be considered in the numerical test section.

Let us consider the orthogonal basis 
\begin{equation}
	\label{eq:bases_asim}
	\{\phi_0,\phi_0,\dots,\phi_N\} \quad\mbox{and}\quad \{\psi_0,\psi_0,\dots,\psi_M\} \qquad\qquad\mbox{with}\quad  \phi_0=\psi_0=1
\end{equation}$\{\phi_j\}_{j=0}^N$ and  $\{\psi_i\}_{i=0}^M$, with $\phi_0=\psi_0=1$, and the the velocity and concentration expansions
$$\tilde{u}(t,x,\xi) = \sum_{j=0}^{N} u_j(t,x) \phi_j(\xi)\qquad\mbox{and}\qquad\tilde{c}(t,x,\xi) = \sum_{i=0}^{M} c_i(t,x) \psi_i(\xi).$$

\begin{table}[!ht]
	\centering
	\begin{tabular}{lccr}
		\hline
		\textbf{Coefficient} & \textbf{Integral Definition} & & \textbf{Dimension} \\ \hline
		$\mu_{k}^u$ & $\int_0^1 \phi_k \phi_k \, d\xi$ & & $(N+1)$ \\
		$\mu_{k}^c$ & $\int_0^1 \psi_k \psi_k \, d\xi$ & & $(M+1)$ \\
		$A^u_{ijk}$ & $\int_0^1 \phi_i \phi_j \phi_k \, d\xi$ & & $(N+1) \times (N+1) \times (N+1)$ \\
		$A^c_{ijk}$ & $\int_0^1 \psi_i \phi_j \psi_k \, d\xi$ & & $(M+1) \times (N+1) \times (M+1)$ \\
		$B^u_{kji}$ & $\int_0^1 \phi_k' \phi_i \left( \int_0^\xi \phi_j \, ds \right) d\xi$ &  &$(N+1) \times (N+1) \times (N+1)$ \\
		$B^c_{kji}$ & $\int_0^1 \psi_k' \psi_i \left( \int_0^\xi \phi_j \, ds \right) d\xi$ & & $(M+1) \times (N+1) \times (M+1)$ \\
		$D_{ki}$ & $\int_0^1 \phi_k \left( \int_\xi^1 \psi_i \, d\tau \right) d\xi$ & & $(N+1) \times (M+1)$ \\
		$E_{ki}$ & $\int_0^1 \phi_k \xi \psi_i \, d\xi$ & & $(N+1) \times (M+1)$\\
		\hline
	\end{tabular}
	\caption{\label{tab:coefAsimetrico}System coefficients for the asymmetric moment system.}
\end{table}

Defining the coefficients in Table \ref{tab:coefAsimetrico}, where many of them are no longer symmetric since mixed integrals between the families $\{\phi\}$ (velocity) and $\{\psi\}$ (sediments) appears, the final shallow moment system with respect to basis \eqref{eq:bases_asim}  reads, for  $k=0,1,\dots,N$ and $l=0,1,\dots,M$,
	\begin{empheq}[left={\empheqlbrace\,}]{align}
		&\partial_t h + \partial_x (hu_0) = 0,\nonumber\\[2mm]
		& \mu^u_k\partial_t (h u_k) + \partial_x \left( h \sum_{i,j=0}^N u_i u_j A_{ijk}^u \right) + \sum_{i,j=0}^N u_i \partial_x (h u_j) B^u_{kji} - \partial_x (h u_0) \sum_{i=0}^N u_i B^u_{k0i} \nonumber\\ 
		&\qquad+gh \partial_x (b+h) \delta_{k0} + g \Delta \rho  \left[h\dsum_{i=0}^{M} D_{ki}\p_x(hc_i) + \partial_x h  \sum_{i=0}^M E_{ki}\, (hc_i) + \partial_x b  \sum_{i=0}^M A^c_{0ki}\, (hc_i)\right] = 0, \nonumber
		\\[2mm]
		& \mu_l^c\partial_t (h c_{l}) + \partial_x \left(h \sum_{i=0}^M\sum_{j=0}^N c_i u_j A^c_{ijl}\right) + \sum_{i=0}^M\sum_{j=0}^N c_i \partial_x (h u_j) B^c_{lji} - \partial_x(h u_0) \sum_{i=0}^M c_i B^c_{l0i} + \mathcal S^c_l(c)
		=0.\nonumber
	\end{empheq}
with
$$
\mathcal S_l(\mathbf c) =
\int_0^1 \psi_l(\xi) \partial_\xi F(\tilde c(\xi))\,d\xi .
$$

An interesting case will be the use of a unique basis $	\{\phi_0,\phi_0,\dots,\phi_M\}$ for both variables (velocity and concentration), but considering different subspaces to project each variable. Let us consider $N\leq M$ and the sediment concentration expanded in the full space $V_M$ of dimension $M+1$, while the velocity is restricted to the subspace $V_N \subseteq V_M$ of dimension $N+1$. In this case, the resulting system is simpler. It reads for $k=1,...,N$, $l=1,...,M$
\begin{empheq}[left={\empheqlbrace\,}]{align}
	&\partial_t h + \partial_x (hu_0) = 0,\nonumber\\[2mm]
	& \mu_k\partial_t (h u_k) + \partial_x \left( h \sum_{i,j=0}^N u_i u_j A_{ijk} \right) + \sum_{i,j=0}^N u_i \partial_x (h u_j) B_{kji} - \partial_x (h u_0) \sum_{i=0}^N u_i B_{k0i} \label{eq:sistemNoUnif}\\ 
	&\qquad+gh \partial_x (b+h) \delta_{k0} + g \Delta \rho  \left[h\dsum_{i=0}^{M} D_{ki}\p_x(hc_i) + \partial_x h  \sum_{i=0}^M E_{ki}\, (hc_i) + \partial_x b  \sum_{i=0}^M A_{ki0}\, (hc_i)\right] = 0, \nonumber
	\\[2mm]
	& \mu_l\partial_t (h c_{l}) + \partial_x \left(h \sum_{i=0}^M\sum_{j=0}^N c_i u_j A_{ijl}\right) + \sum_{i=0}^M\sum_{j=0}^N c_i \partial_x (h u_j) B_{lji} - \partial_x(h u_0) \sum_{i=0}^M c_i B_{l0i} + \mathcal S_l(c)
	=0.\nonumber
\end{empheq}
with
$$
\mathcal S_l(\mathbf c) =
\int_0^1 \phi_l(\xi) \partial_\xi F(\tilde c(\xi))\,d\xi .
$$

The advantage of this formulation is that, since a common basis is used, the coefficients need to be computed only once over the maximum index range $[0,M]$. Moreover, this approach provides additional flexibility in the formulation of the system. In particular, when several variables are approximated, different numbers of modes may be used whenever a more accurate representation of one variable is required.
In our case, it is expected that more function bases are needed to reproduce the vertical profile of concentrations, which include constant areas, whereas the fluid could be good enough represented by only a few bases. It will be illustrated in the numerical test section (see Section \ref{se:numericalTest}). 

\subsection{A particular basis: $\mathbb{P}_1$ functions}\label{se:ParticularBasis}

Legendre polynomials have been widely used as a polynomial basis to derive shallow moment models (see, e.g., \cite{Kowalski2018, koellermeier2020analysis, Koellermeier2020,GarresDiaz2021b}). This polynomial bases are highly accurate for smooth configurations, such as the parabolic velocity profiles of viscous fluids. However, in many applications, the observed vertical distributions (of velocity, sediment concentration, or other variables) exhibit sharp gradients or nearly constant zones. Standard polynomials fails to reproduce these features efficiently due to their inherently oscillatory nature (see \cite{GarresDiaz2023}). Alternatively, in \cite{Koellermeier2025}, first- and second-order splines were introduced to reconstruct the vertical profile of velocity. Nevertheless, because the resulting basis was non-orthogonal, a coefficient matrix had to be inverted, which significantly increased the complexity of the final system, concretely within the non-conservative product terms.

In the previous subsection, the shallow moment model was derived for an arbitrary orthogonal basis $\{\phi_0, \phi_1, \dots, \phi_N\}$. In this section, with the purpose of accurately capturing the sharp vertical profiles of sediment concentration, we introduce a specific basis capable of significantly improving these representations compared to previous moment models in the literature. 

We consider a functional basis inspired by the classical continuous piecewise linear polynomials ($\mathbb{P}_1$) which one usually finds in the finite element framework. In order to gain flexibility in our approach and allows for a local grid refinement, we establish a general, possibly non-uniform, partition of the vertical domain $\xi \in [0,1]$ into $N$ subintervals. Let $0 = \xi_0 < \xi_1 < \dots < \xi_N = 1$ denote the grid nodes, with grid spacings defined as $\Delta \xi_{i+1/2} = \xi_{i+1} - \xi_i$ for $i=0, \dots, N-1$. Associated with this partition, we introduce a set of nodal functions $\{\varphi_i(\xi)\}_{i=0}^N$ satisfying $\varphi_i(\xi_j) = \delta_{ij}$. Specifically, for each internal node $i = 1, \dots, N-1$, these basis functions take the form
\begin{equation}\label{eq:P1_basis}
	\varphi_i(\xi) = 
	\begin{cases} 
		\dfrac{\xi - \xi_{i-1}}{\Delta \xi_{i-1/2}}, & \text{if } \xi \in [\xi_{i-1}, \xi_i], \\
		\dfrac{\xi_{i+1} - \xi}{\Delta \xi_{i+1/2}}, & \text{if } \xi \in [\xi_i, \xi_{i+1}], \\
		0, & \text{otherwise}.
	\end{cases}
\end{equation}
and 
\begin{equation}
	\varphi_0(\xi) = 1 
	\quad \text{and} \quad
	\varphi_N(\xi) = 
	\begin{cases} 
		\dfrac{\xi - \xi_{N-1}}{\Delta \xi_{N-1/2}}, & \text{if } \xi \in [\xi_{N-1}, \xi_N], \\
		0, & \text{otherwise}.
	\end{cases}
\end{equation}

Finally, we apply the Gram-Schmidt orthogonalisation procedure numerically as a pre-processing step to the set $\{\varphi_k\}_{k=0}^N$. Since the first function is the constant mode, the resulting basis satisfies $\phi_0(\xi)=1$. The resulting basis functions $\{\phi_k\}_{k=0}^N$ are continuous, piecewise linear, and orthogonal. Consequently, all the computations derived in the previous subsection remain valid. This strategy avoids the complexity of non-orthogonal systems while gaining the targeted physical flexibility to cluster nodes near critical regions, such as the bottom boundary layer. We shall remark that this orthogonalisation process is performed entirely in a numerical manner as a one-time pre-processing step. All the resulting coefficients are pre-computed, stored, and subsequently read during the simulation without any need for recalculation, therefore, no additional computational efforts are needed.

\section{Numerical approximation}\label{se:numericalApprox}
In order to describe the numerical approximation of system \eqref{eq:finalModel}, it is convenient to write the system in a compact form. Let us consider $$\bd{U} = \left(U_0,...,U_N\right)^T\in\mathbb{R}^{N+1},\quad \bd{C} = \left(C_0,...,C_N\right)^T\in\mathbb{R}^{N+1},\qquad\mbox{with}\quad U_k = hu_k,\ C_k=hc_k,$$ for $k=0,\dots,N$. We define now the vector of conservative variables  $\bd{W} = \left(h, \bd{U}^T, \bd{C}^T, b\right)^T \in\mathbb{R}^{2N+4}$. Then, system \eqref{eq:finalModel} is written in a compact form as
\begin{equation}
	\label{eq:SWEM}
	\p_t \bd{W} + \p_x \bd{F}(\bd{W}) + \bd{B}(\bd{W})\p_x \bd{W} + \bd{E}(\bd{W}) =0,
\end{equation}
where $\bd{F}$ is the conservative flux, $\bd{B}$ is the matrix containing non-conservative products, and $\bd{E}$ is the source terms representing sedimentation terms. In next, we detail these terms. 

The conservative flux is defined by $\bd{F}(\bd{W})=\left(U_0, \bd{F}_U(\bd{W})^T, \bd{F}_C(\bd{W})^T,0\right)^T$ with $\bd{F}_U,\bd{F}_C\in\mathbb{R}^{N+1}$, whose $k$-th component is defined by
\begin{equation*}
	(\bd{F}_{U})_{k+1} = \dfrac{1}{\mu_k}\left(g\dfrac{h^2}{2}\delta_{k0} + \dfrac{1}{h} \dsum_{i,j=0}^N U_i U_j A_{ijk}\right),\qquad (\bd{F}_{C})_{k+1} = \dfrac{1}{\mu_k\, h} \dsum_{i,j=0}^N C_i U_j A_{ijk},\quad\mbox{for }k=0,\dots, N.
\end{equation*}
The non-conservative product term is $\bd{B}(\bd{W})\p_x \bd{W}$ defined through the matrix $\bd{B}(\bd{W}) \in \mathbb{R}^{(2N+4)\times(2N+4)}$, that is written in block-wise form as
\begin{equation*}
	\bd{B}(\bd{W}) = 
	\begin{pmatrix}
		0 & \bd{0}^T & \bd{0}^T & 0 \\[2mm]
		\bd{B}_{U,h} & \bd{B}_{U,U} & \bd{B}_{U,C} & \bd{B}_{U,b} \\[2mm]
		\bd{0} & \bd{B}_{C,U} & \bd{O} & \bd{0} \\[2mm]
		0 & \bd{0}^T & \bd{0}^T & 0
	\end{pmatrix},
\end{equation*}
with $\bd{O}\in\mathbb{R}^{N+1}\times\mathbb{R}^{N+1}$ and $\bd{0}\in\mathbb{R}^{N+1}$ are the null matrix and vector, respectively. The non-zero vector blocks $\bd{B}_{U,h}, \bd{B}_{U,b} \in \mathbb{R}^{N+1}$ and $\bd{B}_{U,U}, \bd{B}_{U,C}, \bd{B}_{C,U} \in \mathbb{R}^{(N+1) \times (N+1)}$ are defined element-wise for indices $k, j = 0, \dots, N$ as:
\begin{equation*}
	(\bd{B}_{U,h})_{k+1} =  \frac{g \Delta \rho}{\mu_k} \sum_{i=0}^N E_{ki}\, C_i,\qquad (\bd{B}_{U,U})_{k+1,j+1} = \frac{1}{\mu_k\,h} \sum_{i=0}^N  \left(B_{kji} - \delta_{j0}B_{k0i}\right)U_i,
\end{equation*}
		\begin{equation*}
			(\bd{B}_{U,C})_{k+1,j+1} = \frac{g \Delta \rho h}{\mu_k} D_{kj},\qquad
	(\bd{B}_{U,h})_{k+1} =  \frac{g \Delta \rho}{\mu_k} \sum_{i=0}^N E_{ki}\, C_i,\qquad (\bd{B}_{U,b})_{k+1} = \frac{gh}{\mu_k} \delta_{k0} + g \Delta \rho\, C_k.
\end{equation*}
and
\begin{equation*}
	(\bd{B}_{C,U})_{k+1,j+1} = \frac{1}{\mu_k\,h}  \sum_{i=0}^N  \left(B_{kji} - \delta_{j0}B_{k0i}\right)C_i.
\end{equation*}
and the source term described by $\bd{E}\in\mathbb{R}^{2N+4}$ whose non-zero components are 
\begin{equation*}
	(\bd{E})_{N+3+k} = \dfrac{\mathcal{S}_k(c)}{\mu_k}, \qquad \mbox{for }k=0,...,N.
\end{equation*}

\medskip 

For the numerical resolution, the proposed system is discretized within the finite volume framework employing a two-step splitting strategy. In the first step, the vertical sedimentation terms are neglected, and the resulting subsystem is discretized and solved within the framework of hyperbolic systems with non-conservative products; in the second step, the sedimentation dynamics are resolved and updated by discretizing the sedimentation terms through a vertical reconstruction--evolution--projection algorithm.

Before describing these steps, some comments about the hyperbolicity of Shallow Water Moment systems should be said. It is well known that shallow moment models are not hyperbolic for $N>1$. Then, model \eqref{eq:finalModel} is not expected to be hyperbolic in general. Some \textit{hyperbolic moment} models have been introduced previously in the literature (see e.g. \cite{Koellermeier2022}) , which are based on removing high-order terms in the matrix system. Although these versions allows us to exactly compute the system eigenvalues, they are still controversial since the system and its solutions are modified. They are only justified in certain regimes where high-order moments contributions are negligible. Nevertheless, some of these versions might be useful to get approximations of the eigenvalues that one can use in practice. In particular, we take into account the eigenvalues approximation of the \textit{Linearized Moment} system (see \cite{Koellermeier2022}), where the lowest ($\lambda^{\rm{Lin}}_-$) and largest ($\lambda^{\rm{Lin}}_+$) eigenvalues are
\begin{equation}
	\label{eq:eigen_linearized}
\lambda_{\pm}^{\rm{Lin}} =u_0 \pm \sqrt{gh + 3\sum_{i=k}^{N}\mu_k u_k^2}.
\end{equation}
Notice that, for any orthogonal functional basis, the correction to the hydrostatic wave speed is determined by the velocity deviation from the mean velocity
\[
\sum_{k=1}^{N}\mu_k u_k^2=\int_0^1\left(\widetilde u(\xi)-u_0\right)^2\,d\xi.
\]

Here, we prefer to use the full shallow moment model even if it is not globally hyperbolic. We refer the reader to \cite{koellermeier2020analysis} for a detailed studied of this fact. However, no numerical issues related to the lack of hyperbolicity have been found in the numerical simulations performed. In particular, when solving numerically the system, we use a HLL type solver which only needs upper and lower bounds of the eigenvalues, contrary to a Roe-type solver. 

\subsection{Step 1: Hyperbolic system with non-conservative products} The system resulting of neglecting sedimentation terms ($\bd{E}(\bd{W})$) can be written in quasi-linear form as
\begin{equation}\label{eq:compact_form}
	\p_t \bd{W} + \bd{A}(\bd{W})\p_x \bd{W} = 0, \qquad\mbox{with}\quad \bd{A}(\bd{W}) = \dfrac{\partial\bd{F}(\bd{W})}{\partial\bd{W}} + \bd{B}(\bd{W}),
\end{equation}
and the path-conservative framework for finite volume methods can be applied to design the numerical scheme approximating the system (see e.g. \cite{Pares2006,Castro2017}). To this aim, the horizontal domain is discretized in control volumes $V_i\in[x_{i-1/2},x_{i+1/2}]$, for $i\in\mathcal{I}$, with uniform mesh-step $\Delta x$. At time $t^n$, the cell average is defined by 
\begin{equation*}
	\bd{W}^n_i = \dfrac{1}{\Delta x}\dint^{x_{i+1/2}}_{x_{i-1/2}} \bd{W}(t^n,x) dx.
\end{equation*} 
and the discretization of system \eqref{eq:SWEM} is 
\begin{equation} \label{eq:step1}
	\begin{array}{l}
		\bd{W}_{i}^{n+1/2}=\bd{W}_i^n +\dfrac{\Delta t}{\Delta x} \left( \bd{\mathcal{F}}_{i-1/2}^n-\bd{\mathcal{F}}_{i+1/2}^n \,
		+ \dfrac12\left(\bd{\mathcal{B}}_{i+1/2}^n+\bd{\mathcal{B}}_{i-1/2}^n \right)\right),
	\end{array}\end{equation}
with 
\begin{equation*}
	\begin{array}{l}
		\bd{\mathcal{F}}_{i+1/2} = \dfrac{1}{2} \left(\bd{F} (\bd{W}_i)+\bd{F}(\bd{W}_{i+1})\right) \,
		- \,\dfrac{1}{2}\bd{\mathcal{D}}_{i+1/2} \\
	\end{array}
\end{equation*}
and $\bd{\mathcal{B}}_{i+1/2}$ consistently defined in the path-conservative framework. Concretely, 
\begin{equation*}
		\bd{\mathcal{B}}_{i+1/2} =\dfrac12 \dsum_{q=1}^2 \bd{B} \left(\bd{W}_{i+1/2}^{(q)}\right)\left(\bd{W_{i+1}}-\bd{W_{i}}\right),
	\end{equation*}
corresponding to the approximation using the two-point Gauss--Legendre quadrature rule along the straight-line path connecting $\bd{W}_{i}$ and $\bd{W}_{i+1}$. The corresponding quadrature states are
\[
\bd{W}_{i+1/2}^{(q)}
=
(1-s_q)\bd{W}_{i}+s_q\bd{W}_{i+1},
\qquad
s_{1,2}=\frac12\mp\frac{\sqrt{3}}{6}.
\]
Finally, $\bd{\mathcal{D}}_{i+1/2}$ is the numerical diffusion of the scheme. Here we use a HLL scheme, which is defined as a PVM scheme (see \cite{CastroDiaz2012}) by
\begin{equation*}\label{eq:diffusion2}
	\bd{\mathcal{D}}_{i+1/2} = \alpha^0_{i+1/2} \left(\bd{W}_{i+1}-\bd{W}_i -\bd{S}^*_{i+1/2} \right)+ \alpha^1_{i+1/2}\left(\bd{F}(\bd{W}_{i+1})-\bd{F}(\bd{W}_{i}) + \bd{\mathcal{B}}_{i+1/2}\right), 
\end{equation*}
being $\bd{S}^*_{i+1/2}=(b_{i+1}-b_i,0,\dots,0)$ and
\[
\alpha^0_{i+1/2} = \dfrac{S_R|S_L|-S_L|S_R|}{S_R-S_L},
\hspace{0.2cm}
\alpha^1_{i+1/2} = \dfrac{|S_R|-|S_L|}{S_R-S_L},
\]
being $S_L$ and $S_R$ approximations of the minimum and maximum wave speed. At the interface $x_{i+1/2}$, they are defined by
\[
(S_L)_{i+1/2}
=
\frac12\sum_{q=1}^{2}
\lambda_-\left(\bd{W}_{i+1/2}^{(q)}\right),
\qquad
(S_R)_{i+1/2}
=
\frac12\sum_{q=1}^{2}
\lambda_+\left(\bd{W}_{i+1/2}^{(q)}\right).
\]
Here, $\lambda_{\pm}$ are the approximated system eigenvalues that are detailed below (see \eqref{eq:autovalores}).\\

\noindent\textbf{Eigenvalues approximation:} As commented above, it is not possible to exactly compute the eigenvalues for the system matrix \eqref{eq:compact_form}. So, we use approximation \eqref{eq:eigen_linearized} but incorporating the effects of the suspended sediment. Thus, we use as eigenvalues approximation for the HLL numerical flux
\begin{equation}
	\label{eq:autovalores}
\lambda_{\pm}=u_0\pm\sqrt{gh(1+\Delta\rho c_0)+3\sum_{i=k}^{N}\mu_k u_k^2}.
\end{equation}

This approximation recovers the relevant exact limits in several cases. If the velocity and concentration are vertically uniform ($N=0$, $u_1=\cdots=u_N=0$, $c_1=\cdots=c_N=0$) it reduces to the exact shallow-water with suspended sediments eigenvalues
\[
u_0\pm\sqrt{gh(1+\Delta\rho c_0)}.
\]
Moreover, \eqref{eq:autovalores} are also an exact approximation of the extremal eigenvalues  for the case of linear velocity $(N=1)$ with a constant sediment profile ($c_1=0$) provided $A_{111}=0$. This condition holds for the shifted Legendre bases and also for the $\mathbb{P}_1$ bases in Subsection \ref{se:ParticularBasis} with $N=1$. For this family of $\mathbb{P}_1$ bases, in the uniform case, one gets
\[
A_{111}=\frac{N-1}{4N^2},
\]
so that \(A_{111}=0\) for \(N=1\) and $A_{111}\longrightarrow0$ as $N\longrightarrow\infty$. That is, it is a small term ($\sim \mathcal{O}(1/N)$) and the approximation used here seems to be reasonable. As mentioned before, no numerical issues related to the possible loss of hyperbolicity have been found in the numerical tests. Thus, the proposed formula \eqref{eq:autovalores} preserves the exact hydrostatic contribution and incorporates the additional deviation from a vertically uniform velocity and the baroclinic contribution to the pressure.

\subsection{Step 2: Vertical sedimentation}\label{se:step2Esquema}
In this step, sedimentation is taken into account in the concentration evolution equation \eqref{eq:finalModel_c}. Then, we trivially have $h^{n+1}_i = h_i^{n+1/2}$ and $\bd{U}^{n+1}_i = \bd{U}_i^{n+1/2}$, and $\bd{C}^{n+1}_i$ must be computed, for $i\in\mathcal{I}$. As has been previously commented on, we consider a particular algorithm.

The motivation for considering a special treatment for this term is the fact that concentration values must be limited to physical admissible values, i.e., $\wtd{c}(\xi)\in [0,c_{\rm{max}}]$ for $\xi\in[0,1]$. One could discretize this term in the moment framework as \eqref{eq:sedimentacionMomentos}. However, filtering these concentration values is not easy. Instead, we consider a reconstruction--evolution--projection procedure at each control volume $V_i$ consisting of: (i) the concentration values are projected in a fine vertical mesh; (ii) the vertical transport is computed by an upwind scheme; (iii) the updated concentration profile in the fine mesh is projected onto the moment framework, obtaining the updated values for the moments ($C_k$) for concentration. This procedure is detailed in following lines. 

\medskip

First, a uniform grid $\xi_1<\xi_2<\cdots<\xi_{N_z}$ of $[0,1]$ is considered, where $\xi_j=\left(j-\frac12\right)\Delta \xi$ for $j=1,\dots,N_z$ with $\Delta \xi = \frac{1}{N_z}$. Take into account that this division is performed for the sedimentation step. Although we are using the same notation for simplicity, this is not the same discretization as the one used in \eqref{eq:P1_basis} for the definition of the basis. In practice, we set $N_z=100$. Then, the concentration profile is reconstructed on this grid, at each control volume $V_i$, by defining 
\begin{equation*}
	\wh{c}_j^{\,*}
	=
	\dfrac{1}{h}\sum_{k=0}^{N} C_{k,i}^{\,n+1/2} \phi_k(\xi_j),
	\quad\mbox{and}\quad \wh{c}_j^{\,n+1/2} = \max\{0,\min\{\wh{c}_j^{\,*},c_{\rm{max}}\}\},\qquad
	\mbox{for }j=1,\dots,N_z.
\end{equation*}
Notice that in this step we filter the admissible values of concentration $\wh{c}_j^{\,n+1/2} \in [0,c_{\rm{max}}]$. It is necessary since the reconstructed profile in terms of the moment coefficients may contain non-admissible values of concentration. 

\medskip

The vertical conservation law
\begin{equation*}
	\partial_t \wh{c} + \partial_\xi F(\,\wh{c}\,)=0
\end{equation*}
is solved on the grid by a finite-volume upwind scheme. To this aim, we write 
\begin{equation*}
	\wh{c}_j^{\,n+1}
	=
	\wh{c}_j^{\,n+1/2}
	-
	\frac{\Delta t}{\Delta \xi}
	\left(
	F_{j+1/2}^{n+1/2} - F_{j-1/2}^{n+1/2}
	\right),
\end{equation*}
with $F_{j+1/2}$ the corresponding numerical flux. Since the settling velocity is directed downwards, the numerical flux through the interface $\xi_{j+1/2} = \frac12(\xi_j+\xi_{j+1})$ between vertical cells $j$ and $j+1$ is taken as
\begin{equation*}
	F_{j+1/2}
	=
	-\omega_s \wh{c}_{j+1}
	(1-\wh{c}_j)^m
	\mathbf 1_{\{\wh{c}_j<c_{\max}\}}.
\end{equation*}
It follows the ideas in \cite{Burger2007,MoralesdeLuna2017}, that is, the concentration is transported downwards due to the settling velocity, unless a saturated ares is found below. Therefore, at each interface an upwind flux is used, but taking into account whether the lower cell is saturated or not by means of the saturation factor, which is evaluated in the lower receiving cell. In addition, we consider that there is no exchange of sediment with the bottom and atmosphere, so $F_{1/2}=F_{N_z+1/2}=0$ is fixed. Thus, the sediment is redistributed inside the water column without imposing an explicit deposition flux through the bed.

Note that the parameter $c_{\max}$ is used as a numerical admissibility threshold. More precisely, it prevents sediment from entering a cell whose concentration has already reached the prescribed maximum admissible value.

Finally, after the vertical finite-volume update, the resulting discrete profile is projected back onto the moment basis. Using the values $\wh{c}_j^{\,n+1}$ computed at the midpoint nodes $\xi_j$, the updated moment coefficients at control volume $V_i$, i.e.,
$\bd{C}_i^{n+1}=(C_{0,i}^{n+1},C_{1,i}^{n+1},\dots,C_{N,i}^{n+1})$, are approximated as
\begin{equation*}
	C_{k,i}^{n+1}
	=
	\frac{h_i^{n+1}}{\mu_k}
	\dint_{0}^1
	\wh{c}^{\,n+1}(\xi)\phi_k(\xi)\,d\xi
	\approx
	h_i^{n+1}\,
	\dfrac{\Delta \xi}{\mu_k}
	\dsum_{j=1}^{N_z}
	\wh{c}_j^{\,n+1}\phi_k(\xi_j),
	\qquad
	\mbox{for }k=0,\dots,N,
\end{equation*}
where the composite midpoint rule is used to approximate the integral in the previous equation.
 This procedure allows the nonlinear vertical settling dynamics to be resolved on a fine vertical grid, while the global moment model remains written in terms of a finite number of vertical moment coefficients. We shall remark that this procedure is carried out within each control volume in a decoupled form, so it can be parallelized. The computational effort corresponding to this step is evaluated in Subsection \ref{subsec:test_bump_sediment}.%

\begin{remark}
	The numerical scheme presented here is stable under the usual CFL stability restriction
	$$
	\Delta t \leq \dfrac{\mbox{CFL}\,\Delta x}{\displaystyle\max_{i\in\mathcal{I}} \lambda_i}, \quad\mbox{with}\quad \lambda_i=\abs{u_0}+\sqrt{gh(1+\Delta\rho c_0)+3\sum_{i=k}^{N}\mu_k u_k^2}
	$$
	and $\mbox{CFL}<1$. In practice we set $\mbox{CFL}=0.9$ for all the numerical tests. We shall remark that the vertical transport detailed in Step 2 should contribute to the stability condition. However, the characteristic speed of the vertical motion, related to the settling velocity, is very small. Therefore, in practice, this contribution is less restrictive than the one of the hyperbolic step, so we do not include it here. Nevertheless, in cases where the settling velocities where large, then the CFL condition should be checked accordingly.
\end{remark}

\section{Numerical experiments}\label{se:numericalTest}

In this section, four test cases are considered to assess the proposed bases. The first two tests (Subsections \ref{subsec:test_bump_clear} and \ref{subsec:test_bump_clear_nounif}) consider clear water with a boundary-layer profile, 
in order to evaluate the ability of the proposed $\mathbb{P}_1$ basis to represent the velocity profile. The remaining two tests (Subsections \ref{subsec:test_sed_vertical} and \ref{subsec:test_bump_sediment}) consider vertical profiles
 of suspended sediment concentration: the first one is horizontally uniform, aimed at showing the ability of the model to reproduce the vertical sedimentation processes;
 the second a more general case, allowing for an enhanced representation of the concentration profile in regions with sharp vertical variations.

The notation used to refer to the moment model with the different vertical bases considered throughout this section as follows:
\begin{itemize}
	\item Leg-N: Polynomial Legendre basis with $N$ functions for velocity and concentration. 
	\item $\mathbb{P}_1$-N,M: $\mathbb{P}_1$ basis with $N$ and $M$ functions for velocity and concentration, respectively. When only a single number $N$ appears ($\mathbb{P}_1$-N), the same number of basis functions is used for velocity and concentration.
	\item NU-$\mathbb{P}_1$-N: $\mathbb{P}_1$ basis with $N$ functions for velocity and concentration, with non-uniform distribution of the nodes.
\end{itemize}

\subsection{$\mathbb{P}_1$-bases and local refinement for boundary-layer profiles}\label{subsec:test_bump_clear}

The purpose of this first numerical test is to compare different vertical bases
for the approximation of a velocity profile exhibiting a thin boundary layer near
the bottom. We first compare the standard Legendre polynomial and the $\mathbb{P}_1$ bases with the same number of basis functions. We then illustrate a key advantage of the $\mathbb{P}_1$ formulation: its ability to incorporate local refinement through a non-uniform partition, a flexibility that is not naturally available with a global polynomial basis.
We consider a flat bottom and a constant initial water height,
\[
b(x)=0,
\qquad
h(x,0)=1.
\]
Using the normalized vertical coordinate $\xi\in[0,1]$, defined by the change of variable \eqref{eq:cambioVariable}, the reference velocity profile is chosen as
\begin{equation}\label{eq:uref_capalimite}
u_{\mathrm{ref}}(\xi)
=
U_\infty
\left[
\left(
2\frac{\xi}{\delta}
-
\left(\frac{\xi}{\delta}\right)^2
\right)\mathbf{1}_{\xi<\delta}
+
\mathbf{1}_{\xi\geq\delta}
\right], \qquad \mbox{with}\qquad U_\infty=0.5,
\quad
\delta=0.15.
\end{equation}
This profile increases rapidly from the bottom and reaches the constant value
\(U_\infty\) above the boundary-layer thickness \(\delta\).

For a given number of velocity moments \(N\), we consider an orthogonal basis $\{\phi_k\}_{k=0}^{N}$ with $\phi_0=1,$ as in the moment formulation introduced above. The projected profile is then
written as
\[
u_{\mathrm{proj},N}(\xi)
=
\sum_{k=0}^{N} \wtd{u}_k \phi_k(\xi).
\]
Here the coefficients \(\wtd{u}_k\) are the moments associated with the
projection of \(u_{\mathrm{ref}}\) onto the chosen approximation space, which are the initial values for the moment coefficients ($u_k(x,0)=\wtd{u}_k$). Since the
purpose of the test is to compare the approximation properties of the different
bases, we show the reconstructed profiles and the projection errors, rather
than the individual modal coefficients.

We begin by comparing the standard shifted Legendre basis with the $\mathbb{P}_1$
basis of the same dimension. In the $\mathbb{P}_1$ basis case, the interval \([0,1]\) is first
divided into \(N\) equal subintervals,
\[
0=\xi_0<\xi_1<\cdots<\xi_N=1,
\qquad
\xi_i=\frac{i}{N},
\qquad i=0,\ldots,N.
\]

\begin{figure}[htbp]
	\centering
		\centering
		\includegraphics[width=0.45\textwidth]{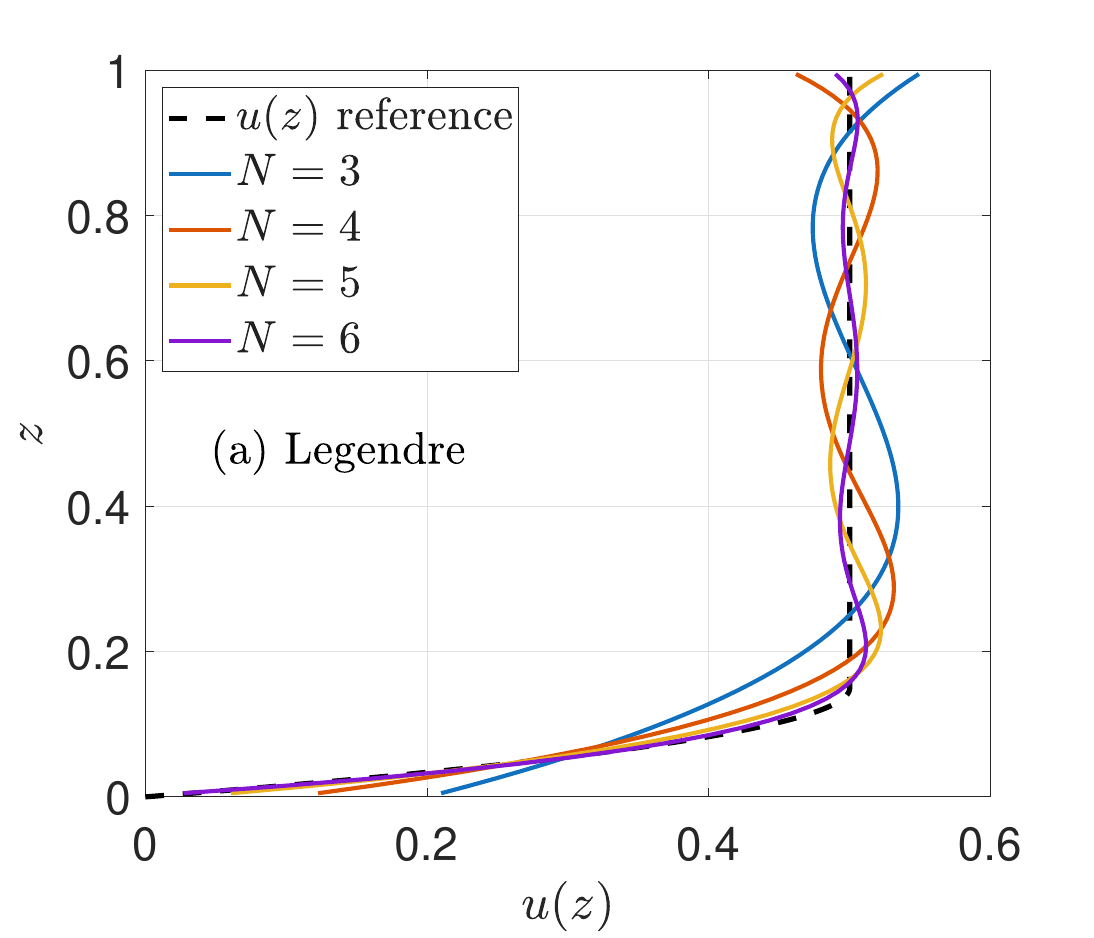}
		\includegraphics[width=0.45\textwidth]{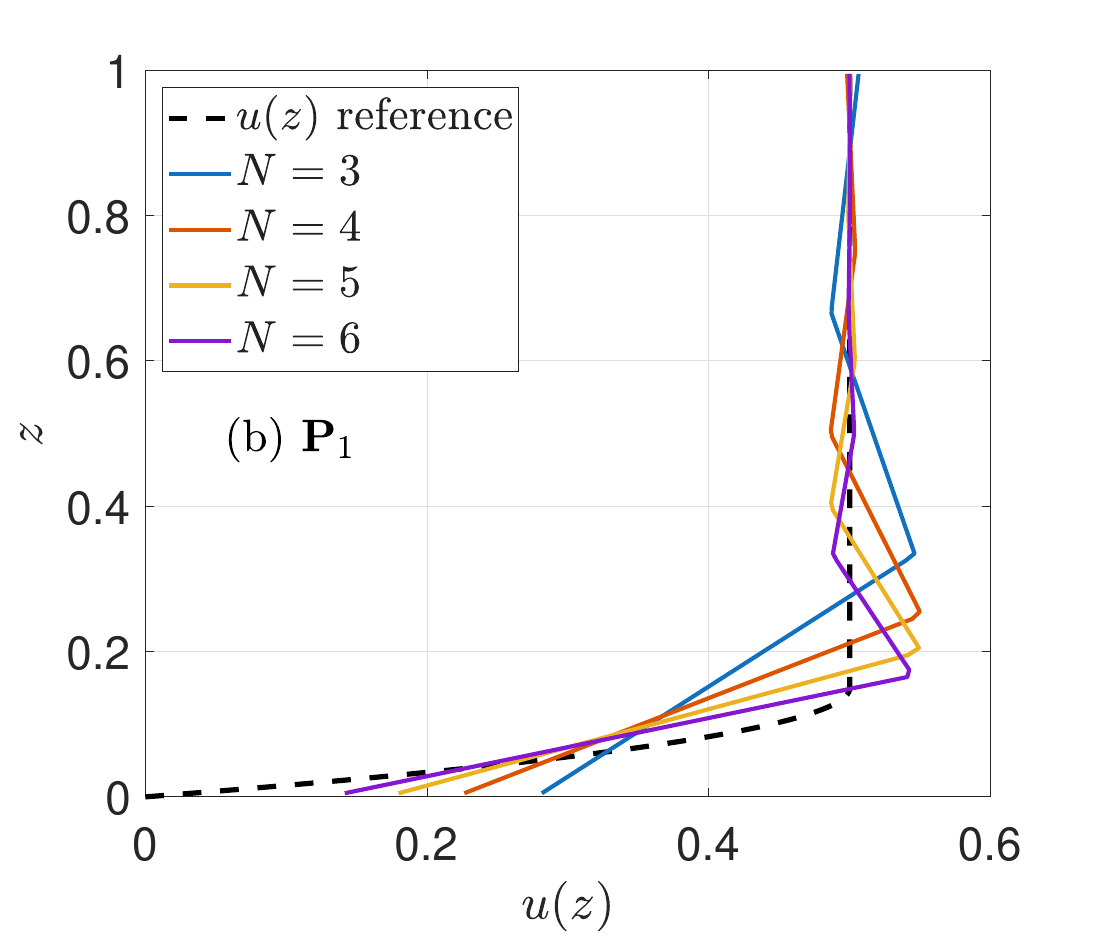}
	\caption{Projection of the boundary-layer profile for \(N=3,4,5,6\) using
		shifted Legendre polynomials (Leg-N, left) and $\mathbb{P}_1$ ($\mathbb{P}_1$-N, right) bases. Both approximations are
		of comparable quality.}
	\label{fig:legendre-vs-uniform-spline}
\end{figure}

Figure~\ref{fig:legendre-vs-uniform-spline} shows the reconstructed profiles for
$N=3,4,5,6$. Qualitatively, both bases families provide comparable
approximations of the boundary-layer profile. Notice that, when the
$\mathbb{P}_1$ basis partition is uniform, neither approximation is specifically adapted to the thin lower layer.

To quantify the approximation, we define
\[
e_N(\xi)
=
u_{\mathrm{proj},N}(\xi)-u_{\mathrm{ref}}(\xi),
\]
and consider the errors in the \(L^\infty\), \(L^2\), and \(L^1\) norms (denoted by $\|e_N\|_{L^\infty}$, $\|e_N\|_{L^2}$, and $\|e_N\|_{L^1}$). %
The corresponding projection errors are shown in
Table~\ref{tab:projection-errors-legendre-uniform-spline}. They confirm that the
Legendre and $\mathbb{P}_1$ approximations are comparable. In particular, for
moderate values of $N$, neither basis offers a clear advantage in this case ($\mathbb{P}_1$ basis constructed on a uniform partition).

\begin{table}[htbp]
	\centering
	\small
	\caption{Projection errors for the boundary-layer profile using shifted Legendre
		and $\mathbb{P}_1$ bases.}
	\label{tab:projection-errors-legendre-uniform-spline}
	\begin{tabular}{c|c|cccc}
		\hline
		Basis & Error & \(N=3\) & \(N=4\) & \(N=5\) & \(N=6\) \\
		\hline
		\multirow{3}{*}{Legendre}
		& \(\|e_N\|_{L^\infty}\)
		& \(1.77\times 10^{-1}\)
		& \(8.99\times 10^{-2}\)
		& \(2.81\times 10^{-2}\)
		& \(1.39\times 10^{-2}\)
		\\
		
		& \(\|e_N\|_{L^2}\)
		& \(4.09\times 10^{-2}\)
		& \(2.52\times 10^{-2}\)
		& \(1.29\times 10^{-2}\)
		& \(5.74\times 10^{-3}\)
		\\
		
		& \(\|e_N\|_{L^1}\)
		& \(2.99\times 10^{-2}\)
		& \(2.00\times 10^{-2}\)
		& \(1.09\times 10^{-2}\)
		& \(4.81\times 10^{-3}\)
		\\
		\hline
		\multirow{3}{*}{$\mathbb{P}_1$}
		& \(\|e_N\|_{L^\infty}\)
		& \(2.49\times 10^{-1}\)
		& \(1.94\times 10^{-1}\)
		& \(1.47\times 10^{-1}\)
		& \(1.09\times 10^{-1}\)
		\\
		& \(\|e_N\|_{L^2}\)
		& \(5.17\times 10^{-2}\)
		& \(4.15\times 10^{-2}\)
		& \(3.20\times 10^{-2}\)
		& \(2.36\times 10^{-2}\)
		\\
		& \(\|e_N\|_{L^1}\)
		& \(3.08\times 10^{-2}\)
		& \(2.33\times 10^{-2}\)
		& \(1.70\times 10^{-2}\)
		& \(1.19\times 10^{-2}\)
		\\
		\hline
	\end{tabular}
\end{table}

The main benefit of the $\mathbb{P}_1$ formulation appears when one takes advantage of
its local character. Unlike global polynomial bases, $\mathbb{P}_1$ bases are associated
with a partition of the interval \([0,1]\). Therefore, the vertical degrees of
freedom can be redistributed according to the expected structure of the profile. To exploit this property, we consider a non-uniform $\mathbb{P}_1$ partition refined near
the bottom. For \(N=3,4,5,6\), the interval \([0,1]\) is divided into \(N\)
subintervals in such a way that the first two subintervals always cover the lower
region \([0,0.1]\), while the remaining \(N-2\) subintervals cover \([0.1,1]\).
Hence,
\[
0=\xi_0<\xi_1<\xi_2=0.1<\xi_3<\cdots<\xi_N=1,
\]
with
\[
\xi_i=\frac{0.1}{2}i,
\qquad i=0,1,2, \qquad \mbox{and}\qquad \xi_{2+j}
=
0.1+j\frac{0.9}{N-2},
\qquad j=1,\ldots,N-2.
\]
The final basis is obtained by orthogonalising the $\mathbb{P}_1$ functions associated with this partition using the Gram--Schmidt procedure. The dimension of the approximation space is still
\(N+1\), but more resolution is now concentrated in the region where the boundary
layer is located.

\begin{figure}[htbp]
	\centering
	\includegraphics[width=0.48\textwidth]{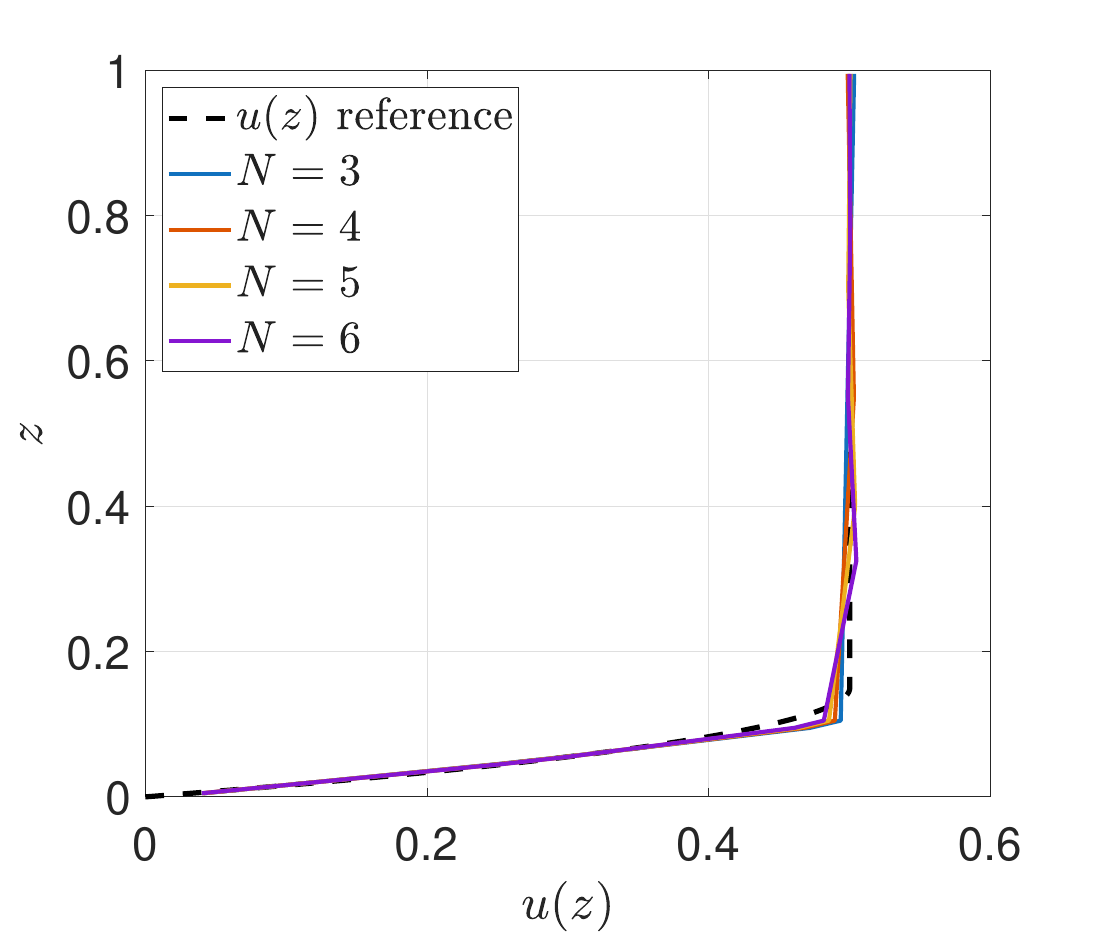}
	\caption{Projection of the boundary-layer profile onto non-uniform $\mathbb{P}_1$
		bases for \(N=3,4,5,6\) (NU-$\mathbb{P}_1$-N). The local refinement near the bottom significantly
		improves the representation of the boundary layer.}
	\label{fig:projection-boundary-layer-spline-nonuniform}
\end{figure}

The resulting projected profiles are shown in
Figure~\ref{fig:projection-boundary-layer-spline-nonuniform}. Qualitatively, the
improvement is clear: the thin bottom layer is much better captured than with both
the Legendre and the $\mathbb{P}_1$ (with uniform partition) bases. The improvement is also confirmed by the projection errors shown in
Table~\ref{tab:projection-errors-nonuniform-spline}. The non-uniform
$\mathbb{P}_1$ bases substantially reduce the errors for all tested values of
$N$. For instance, for $N=3$, the $L^2$ error decreases from
$4.09\times 10^{-2}$ for the Legendre basis and $5.17\times 10^{-2}$ for the
uniform $\mathbb{P}_1$ basis to $7.14\times 10^{-3}$ for the non-uniform
$\mathbb{P}_1$ basis. Thus, in this test case, the decisive improvement does not
come merely from replacing the polynomial basis with a $\mathbb{P}_1$ basis, but
from the possibility of introducing local refinement while retaining the same
number of modes.

\begin{table}[htbp]
	\centering
	\small
	\caption{Projection errors for the boundary-layer profile using non-uniform
		$\mathbb{P}_1$ bases.}
	\label{tab:projection-errors-nonuniform-spline}
	\begin{tabular}{c|cccc}
		\hline
		Error & \(N=3\) & \(N=4\) & \(N=5\) & \(N=6\) \\
		\hline
		\(\|e_N\|_{L^\infty}\)
		& \(3.90\times 10^{-2}\)
		& \(3.52\times 10^{-2}\)
		& \(3.15\times 10^{-2}\)
		& \(2.81\times 10^{-2}\)
		\\
		\(\|e_N\|_{L^2}\)
		& \(7.14\times 10^{-3}\)
		& \(6.68\times 10^{-3}\)
		& \(6.18\times 10^{-3}\)
		& \(5.71\times 10^{-3}\)
		\\
		\(\|e_N\|_{L^1}\)
		& \(3.90\times 10^{-3}\)
		& \(3.72\times 10^{-3}\)
		& \(3.33\times 10^{-3}\)
		& \(2.95\times 10^{-3}\)
		\\
		\hline
	\end{tabular}
\end{table}

This first test therefore motivates the use of $\mathbb{P}_1$ bases in the moment
formulation. If a uniform partition is used, $\mathbb{P}_1$ bases and polynomial bases lead to
comparable approximations. However, $\mathbb{P}_1$ bases offer an additional degree of
flexibility: the vertical partition can be adapted to the expected structure of
the solution. For boundary-layer profiles, this makes it possible to represent the
localized lower layer much more accurately with the same number of moments, an
advantage that is not directly available in the standard global polynomial
setting.

\subsection{Evolution of a boundary-layer velocity profile}\label{subsec:test_bump_clear_nounif}

In the previous test, we compared the ability of several vertical bases to
represent a prescribed boundary-layer velocity profile. We now use the same type of profile as initial condition in the time-dependent moment system and study its evolution, with sediment feedback on the hydrodynamics suppressed by taking $\Delta\rho=0$ in system \eqref{eq:finalModel}. The purpose is to assess whether the improvement obtained by using a
locally refined $\mathbb{P}_1$ basis is also reflected during the dynamical evolution of
the flow.

The computational domain is $x\in[-5,5]$, the bottom is flat and the solution is computed up to $
t_f=1$. The horizontal domain is discretized using \(1000\) cells. The initial water height is given by
\[
h(x,0)
=
1
+
0.5\exp\left(-10(x-1)^2\right)
+
0.5\exp\left(-10(x+1)^2\right).
\]
The initial horizontal velocity is prescribed from the boundary-layer profile \eqref{eq:uref_capalimite}
considered in the previous test with $
U_\infty=1$ and $\delta=0.15$. 

We compare four different velocity approximations: NU-$\mathbb{P}_1$-3, $\mathbb{P}_1$-3, Leg-3 and Leg-6. We recall that, according to the notation, these acronyms denote the non-uniform $\mathbb{P}_1$ and uniform $\mathbb{P}_1$ bases with \(N=3\), and the Legendre bases with \(N=3\) and \(N=6\), respectively. %

The non-uniform $\mathbb{P}_1$ basis is built from the same locally refined partition
used in the previous projection test. Since \(N=3\), the interval \([0,1]\) is
divided into three subintervals, with the first two subintervals concentrated
near the bottom:
\[
0=\xi_0<\xi_1<\xi_2<\xi_3=1,\qquad\mbox{where}\qquad \xi_1=0.05,\ \xi_2=0.1.
\]
The corresponding functions are then orthogonalised, consistently with
the moment formulation used in the model.

\begin{figure}[!h]    \centering\includegraphics[width=0.45\textwidth]{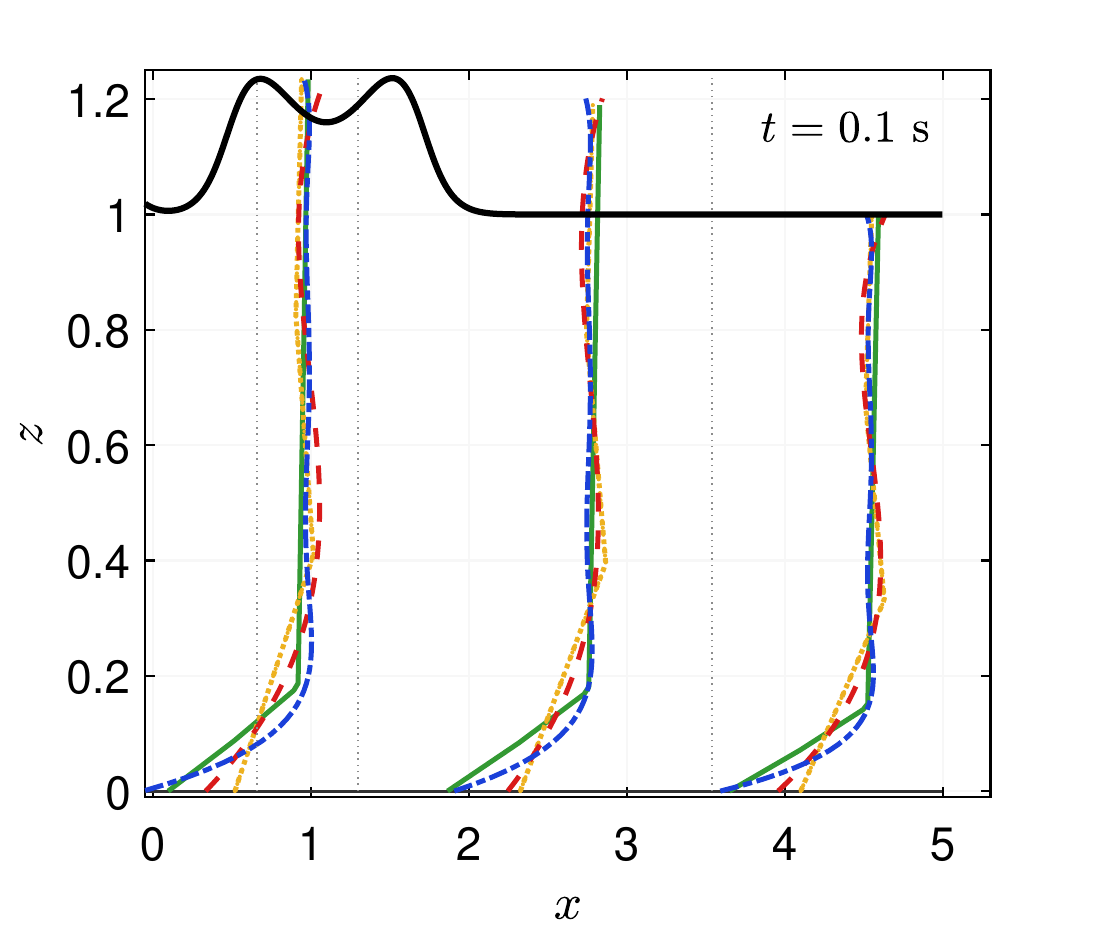}
\includegraphics[width=0.45\textwidth]{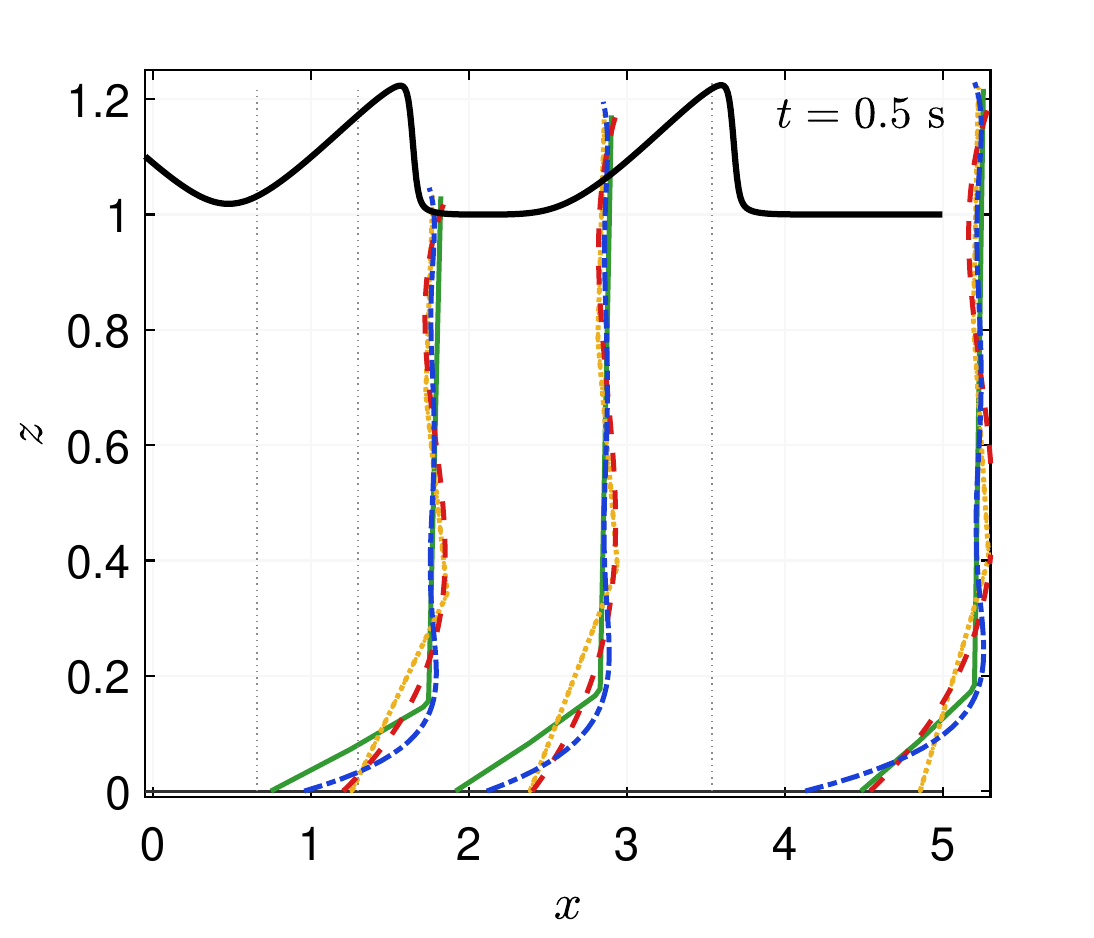}
\includegraphics[width=0.45\textwidth]{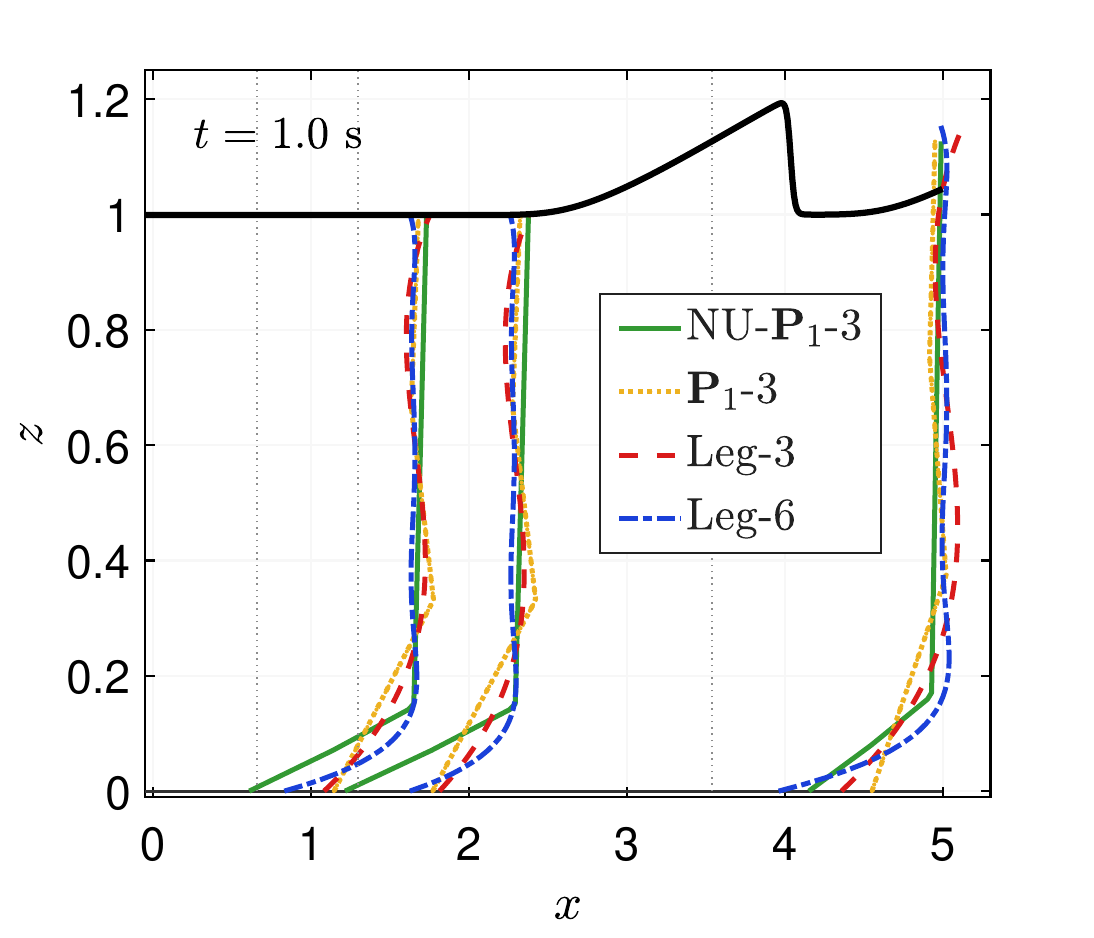}
    \caption{Evolution of the free surface and reconstructed velocity profiles
for the boundary-layer initial condition. The profiles are shown at $t=0.1, 0.5, 1.0$ and the locations \(x=0.655\), \(x=1.295\), and \(x=3.535\) (represented by dotted gray lines).}
    \label{fig:capalim-evolution-full}
\end{figure}

\begin{figure}[!h]    
\includegraphics[width=0.3\textwidth]{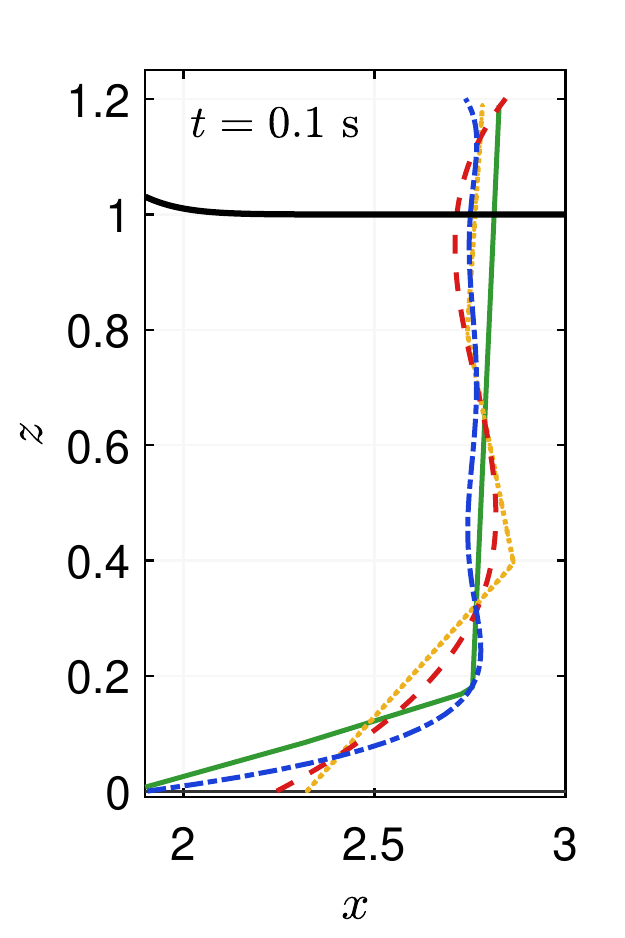}
\includegraphics[width=0.3\textwidth]{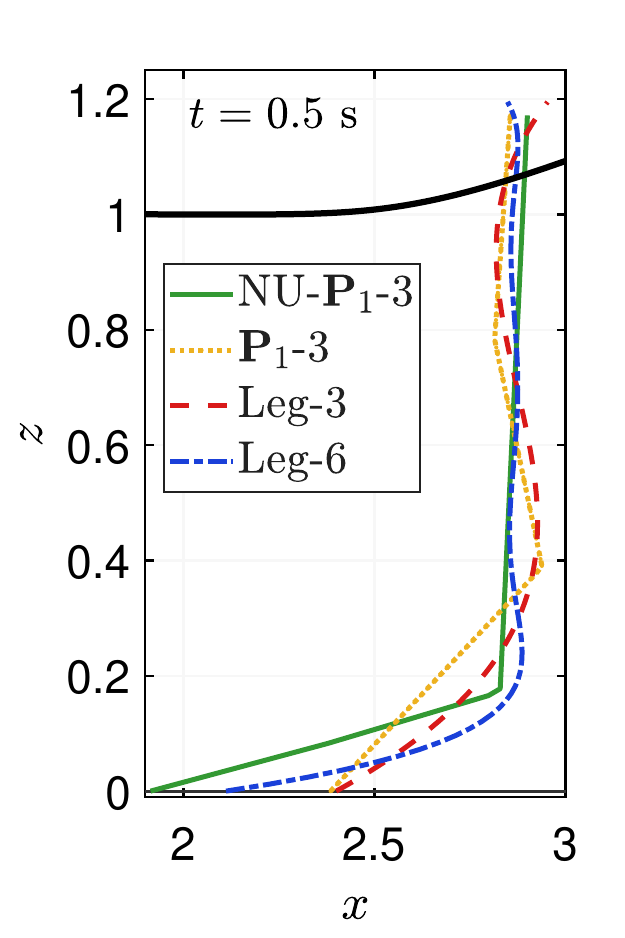}
\includegraphics[width=0.3\textwidth]{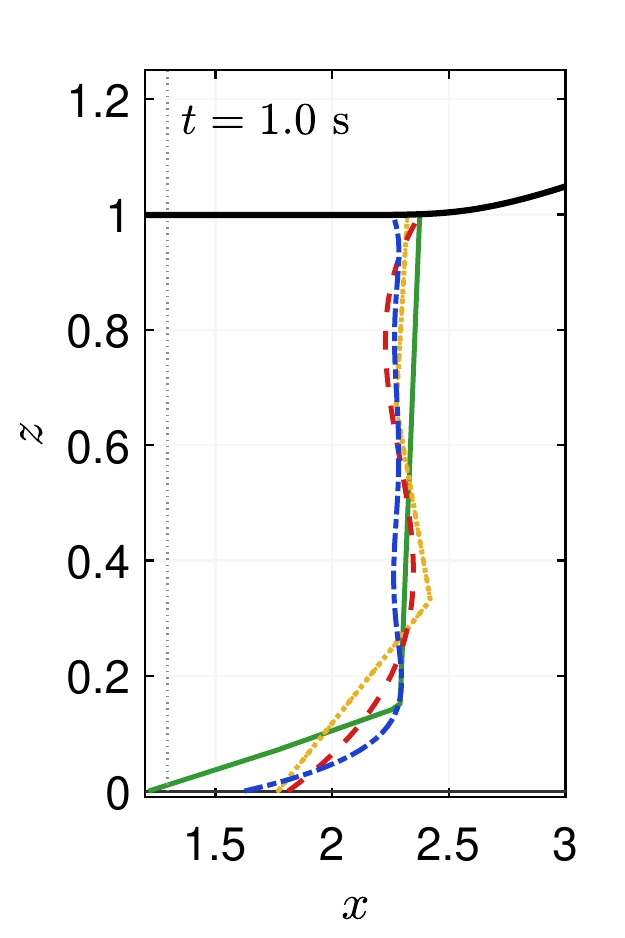}
   \caption{Zoom of the reconstructed velocity profile at \(x=1.295\) and times $t=0.1, 0.5, 1.0$ for the
    boundary-layer evolution test. The non-uniform $\mathbb{P}_1$ basis with \(N=3\)
    gives profiles comparable to those obtained with the Legendre basis with
    \(N=6\).}
    \label{fig:capalim-evolution-zoom}
\end{figure}

In Figures~\ref{fig:capalim-evolution-full} and
\ref{fig:capalim-evolution-zoom}, we show the free surface and the reconstructed
velocity profiles at the times $t=0.1$, $t=0.5$, and $t=1.0$ s. The velocity profiles are displayed at the positions $x=0.655$, $x=1.295$ and $x=3.535$. In the plots, each vertical profile is represented by the curve
\[
x+\tilde u(t,x,z),
\]
so that the horizontal displacement of the curve corresponds directly to the
local horizontal velocity.

Figure~\ref{fig:capalim-evolution-full} shows that the four approximations produce comparable
large-scale dynamics. The main differences are observed near the bottom, where the initial
boundary-layer structure is concentrated. The uniform $\mathbb{P}_1$ basis with \(N=3\)
and the Legendre basis with \(N=3\) provide similar low-order descriptions of the
profile. Both capture the overall deformation of the velocity field, but their
resolution of the lower layer is limited.

The non-uniform $\mathbb{P}_1$ basis with \(N=3\) (NU-$\mathbb{P}_1$-3) gives a sharper representation of
the near-bottom part of the profile. In this region, its behaviour is closer to
that obtained with the higher-order Legendre approximation with \(N=6\) (Leg-6). This shows that
the local refinement near the bottom reproduces part of the effect of adding
global polynomial modes, while keeping only four velocity coefficients. To make this comparison clearer, Figure~\ref{fig:capalim-evolution-zoom} shows a
zoom of the velocity profile around the middle observation region $x=1.295$. We observe that the standard $\mathbb{P}_1$ basis with
\(N=3\) behaves similarly to the Legendre basis with \(N=3\), whereas the
non-uniform $\mathbb{P}_1$ basis with \(N=3\) captures the lower-layer structure in a way
that is comparable to the Legendre approximation with twice as many moments $N=6$. 

\subsection{Vertical sedimentation of an initially uniform concentration}\label{subsec:test_sed_vertical}

We now consider a test designed to isolate the vertical sedimentation mechanism
and to compare the effect of the different bases for the
concentration. The related sedimentation parameters are $c_{\max}=0.8, m=2.5, w_s=0.05.$

The computational domain is $x\in[-5,5]$ and the solution is computed up to $
t_f=100$. The horizontal domain is discretized using \(1000\) cells. Since the dynamics
are purely vertical, all horizontal locations are equivalent and the solution at a single location ($x=0$) will be shown.
The flow is initially at rest, the bottom is flat, and the water
height is constant:
\[
b(x)=0,
\qquad
h(x,0)=1,
\qquad
\tilde u(x,0,\xi)=0.
\]
The initial concentration is assume to be
vertically uniform and given by
\[
\tilde c(x,0,\xi)=c_0,
\qquad
c_0=0.5.
\]

No velocity moments are used in this test. For the modal approximations, we use concentration spaces with \(M=5\). We compare three
different bases: shifted Legendre polynomials (Leg-5), uniform $\mathbb{P}_1$ bases ($\mathbb{P}_1$-5) and non-uniform
$\mathbb{P}_1$ bases (NU-$\mathbb{P}_1$-5). 

The non-uniform $\mathbb{P}_1$ basis is constructed from a partition refined near both
the bottom and the free surface. More precisely, the interval \([0,1]\) is
divided into five subintervals. Two subintervals are placed in the lower region
\([0,0.2]\), one subinterval in the middle region \([0.2,0.8]\), and two
subintervals in the upper region \([0.8,1]\). Thus, the partition is
\[
0=\xi_0<\xi_1<\xi_2<\xi_3<\xi_4<\xi_5=1,
\]
with
\[
\xi_1=0.1,
\qquad
\xi_2=0.2,
\qquad
\xi_3=0.8,
\qquad
\xi_4=0.9.
\]
Finally, the corresponding $\mathbb{P}_1$ functions are orthogonalised.

In addition to the modal solutions, we include a reference solution. This
reference is obtained by solving the vertical sedimentation problem directly on a
fine vertical finite-volume grid in the normalized coordinate
\(\xi\). No modal projection is applied to this solution. Therefore, it provides
a high-resolution approximation of the vertical transport process and serves as
a benchmark for the modal reconstructions. More precisely, if \(c_j^n\) denotes the cell average of the concentration in the
vertical finite-volume grid, the reference is advanced using the same hindered
settling flux as in the vertical sedimentation step (see Subsection \ref{se:step2Esquema}). With zero flux at the bottom
and at the free surface, the update has the form
\[
c_j^{n+1}
=
c_j^n
-
\frac{\Delta t}{\Delta \xi}
\left(
F_{j+1/2}^n-F_{j-1/2}^n
\right),
\]
where
\[
F_{j+1/2}^n
=
-w_s\,c_{j+1}^n
\left(1-c_j^n\right)^m
\mathbf{1}_{c_j^n<c_{\max}}.
\]
The boundary fluxes are set to zero, $F_{1/2}=F_{J+1/2}=0$. 

\begin{figure}[!h]
    \centering
    \includegraphics[width=0.3\textwidth]{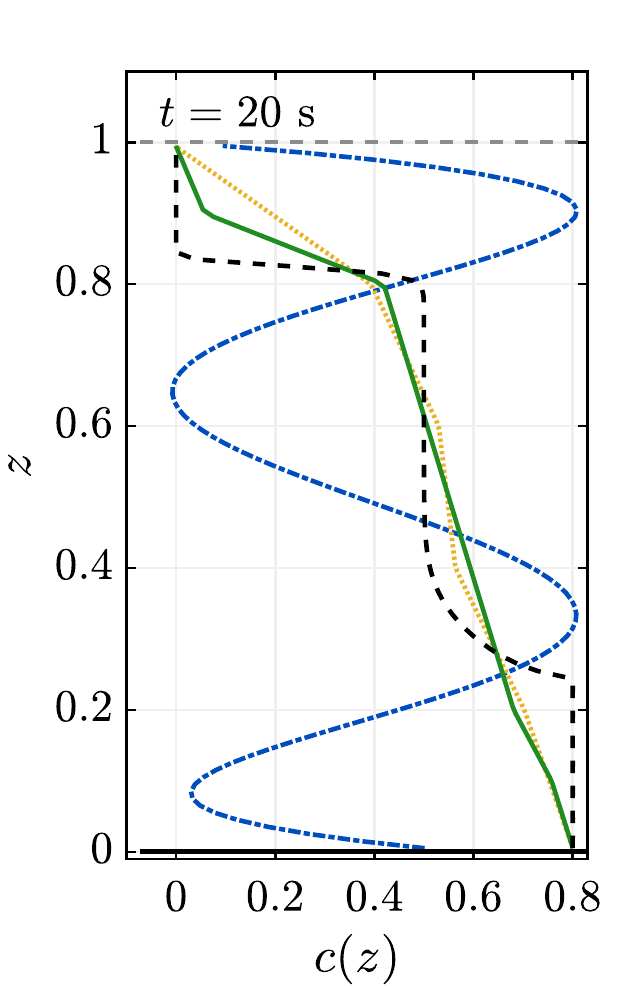}
    \includegraphics[width=0.3\textwidth]{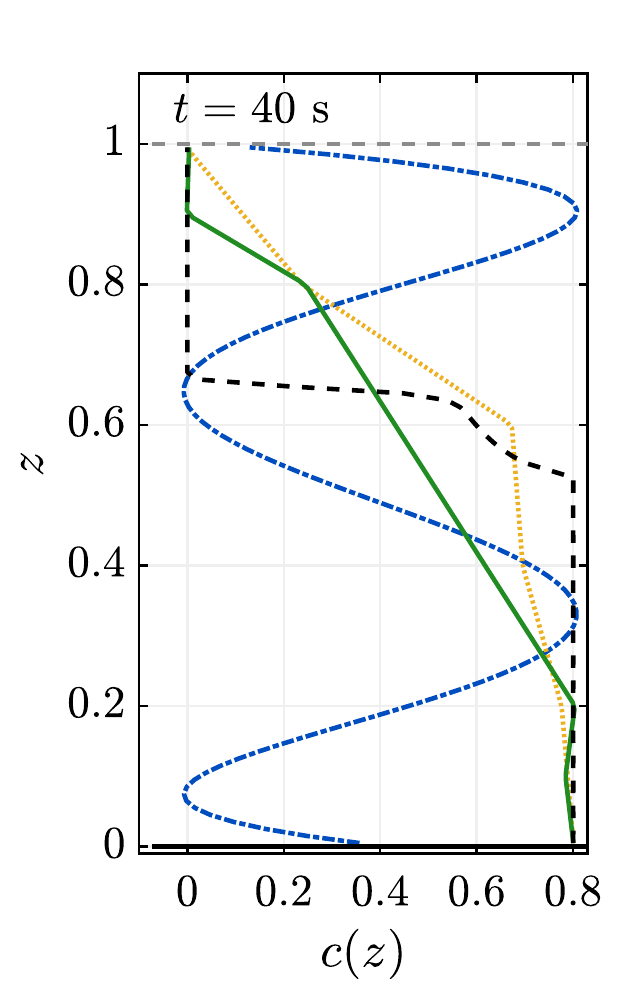}
    \includegraphics[width=0.3\textwidth]{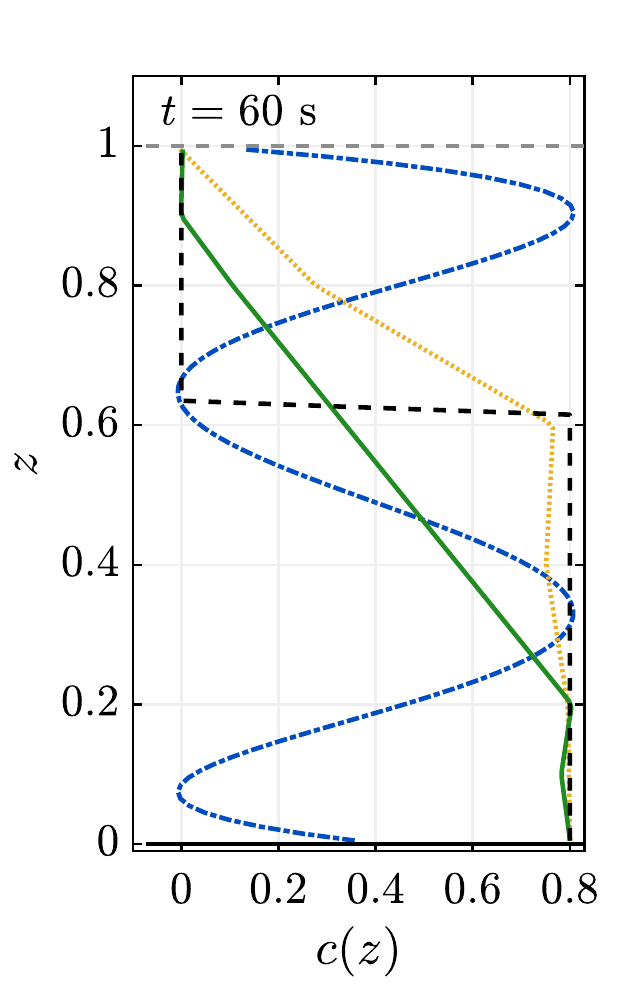}
    \includegraphics[width=0.3\textwidth]{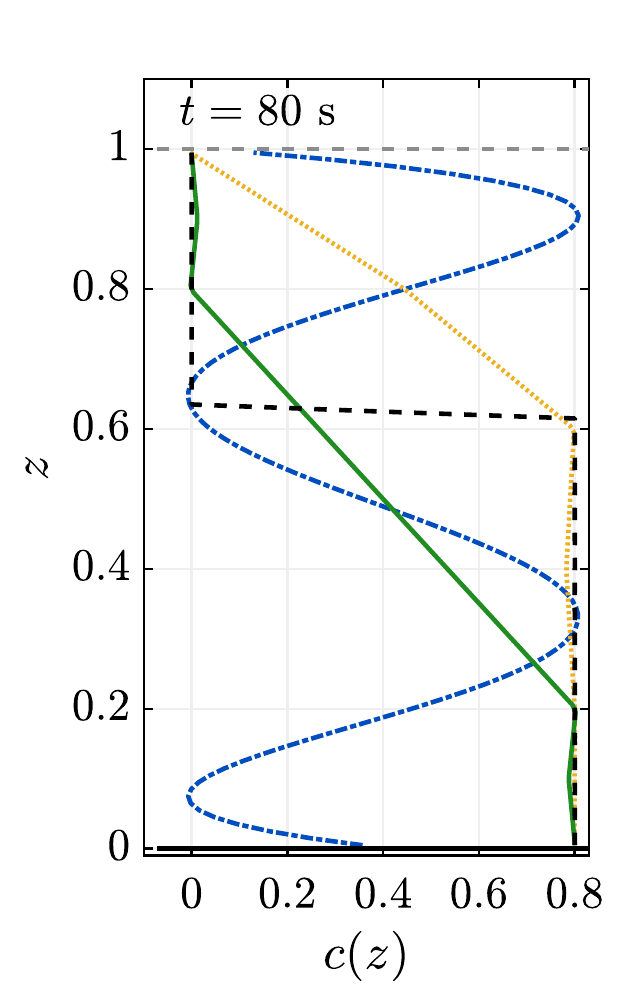}
    \includegraphics[width=0.3\textwidth]{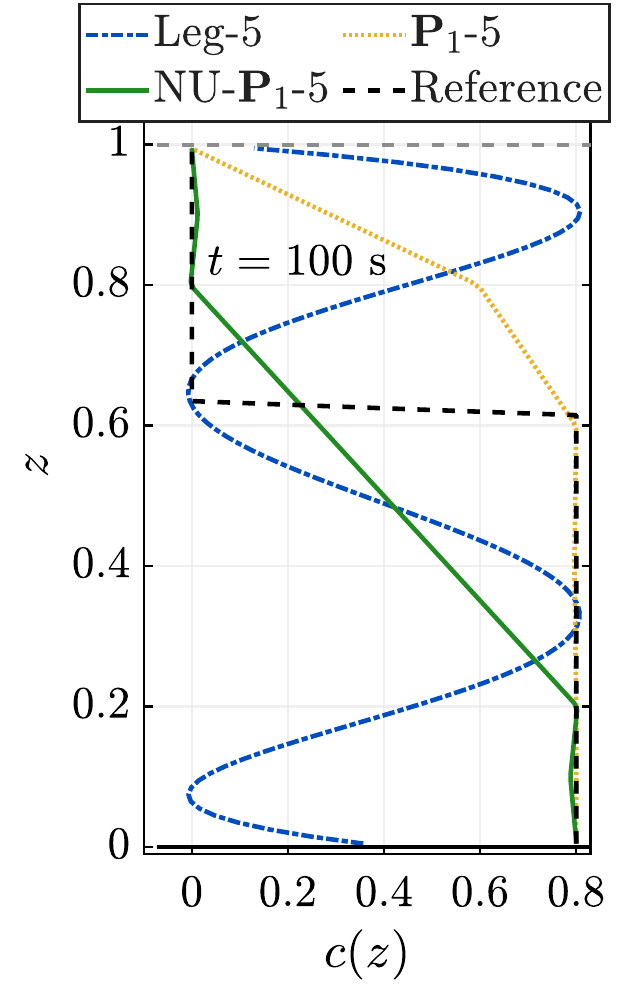}
    \caption{Evolution of the reconstructed concentration profiles for the
    vertical sedimentation test. The initial concentration is uniform,
    \(c_0=0.5\), the maximum concentration is \(c_{\max}=0.8\), and the profiles
    are shown at \(x=0\). The reference solution is computed directly on a fine
    vertical finite-volume grid with \(200\) points, without modal projection.}
    \label{fig:vertical-sedimentation-concentration}
\end{figure}

Figure \ref{fig:vertical-sedimentation-concentration} shows the evolution of the concentration profile at $x=0$ for several times. We see that the reference solution exhibit the expected behaviour of the model. The sediment moves downwards and forms a concentrated lower layer, while the concentration in the upper part of the water column is progressively reduced. At late times, the reference solution approaches a nearly saturated lower layer with concentration close to \(c_{\max}\), together with a diminished upper layer. From mass conservation, the expected height of the saturated region (denoted by $\Delta z_{\textrm{sat}}$) is approximately
\[
\Delta z_{\textrm{sat}}=\frac{h_0 c_0}{c_{\max}}
=
\frac{0.5}{0.8}
=
0.625,
\]
where $h_0$ is the initial fluid depth. This is consistent with the position of the sharp transition observed in the reference profile. 

We also plot in this figure the approximated profiles with the bases Leg-5, $\mathbb{P}_1$ and NU-$\mathbb{P}_1$. The Legendre approximation exhibits strong oscillations as the concentration profile develops sharp vertical gradients. This behaviour, that is far from the physically expected, is typical of global polynomial approximations when the solution contains steep fronts or nearly discontinuous transitions. The uniform $\mathbb{P}_1$ approximation is less oscillatory,
due to the local character of the basis, but the transition region remains
diffused because only a small number of concentration modes is used. The non-uniform $\mathbb{P}_1$ approximation provides a more localized representation of
the vertical structure. In particular, the refinement near the lower and upper
parts of the water column improves the description of the saturated and empty regions (high and low concentration, respectively). 

\subsection{Coupled velocity--sediment transport over a smooth bump}
\label{subsec:test_bump_sediment}

We now consider a test in which both the vertical velocity structure and
the suspended sediment concentration evolve over a non-flat bottom. The
computational domain is $    x\in[-5,5]$ and the bottom topography is given by the smooth bump
\[
    b(x)=0.1\exp\left(-10x^2\right).
\]
The initial free surface is flat, $\eta(x,0)=1$, so that $h(x,0)=1-b(x)$. The computation is performed with 200 cells up to the final time $t_f=100$. In this test we will test the asymmetric formulation presented in Subsection \ref{se:asimetrico}. In particular, we consider the case in system \eqref{eq:sistemNoUnif}, where the velocity field is approximated in a subspace of the functional space used to approximate the sediment concentration. Then, we distinguish between the number of velocity moments and the number of concentration moments. Let us denote by \(N\) the highest velocity mode and by \(M\) the highest concentration mode, with $N\leq M$. For the sediment concentration \(M=5\) is assumed in all simulations. For the velocity field, $N=2$ or $N=5$ will be considered. Concretely, the considered bases are Leg-5 (Legendre basis with $N=M=5$), $\mathbb{P}_1$-5 ($\mathbb{P}_1$ basis with $N=M=5$) and $\mathbb{P}_1$-2,5 ($\mathbb{P}_1$ basis with $N=2$ and $M=5$).

The initial velocity profile is prescribed as the parabolic function
\begin{equation}
	\label{eq:uref}
	    u_{\mathrm{ref}}(\xi)=0.03+0.36\xi-0.18\xi^2,
	\qquad \xi\in[0,1].
\end{equation}
along the whole domain, being the averaged velocity $\int_0^1u_{\mathrm{ref}}(\xi)\,d\xi=0.15.$

We assume that the initial suspended sediment is localized on the left part of the domain while clear water is on the right part. More
precisely,
\begin{equation}
	\label{eq:cref}
c_{\mathrm{ref}}(x,\xi)= \left\{\begin{array}{ll}
	0.05+0.60\xi-0.30\xi^2 & \mbox{if }x < -1,\\
	0 & \mbox{if }x\geq -1, 
\end{array}\right.  \qquad \xi\in[0,1].
\end{equation}
The averaged sediment concentration on the left part of the domain is $\int_0^1 c_{\mathrm{ref}}(\xi)\,d\xi=0.25$. The same sedimentation parameter as in previous test are considered. 

\begin{figure}[!ht]
	\centering
	\includegraphics[width=0.47\textwidth]{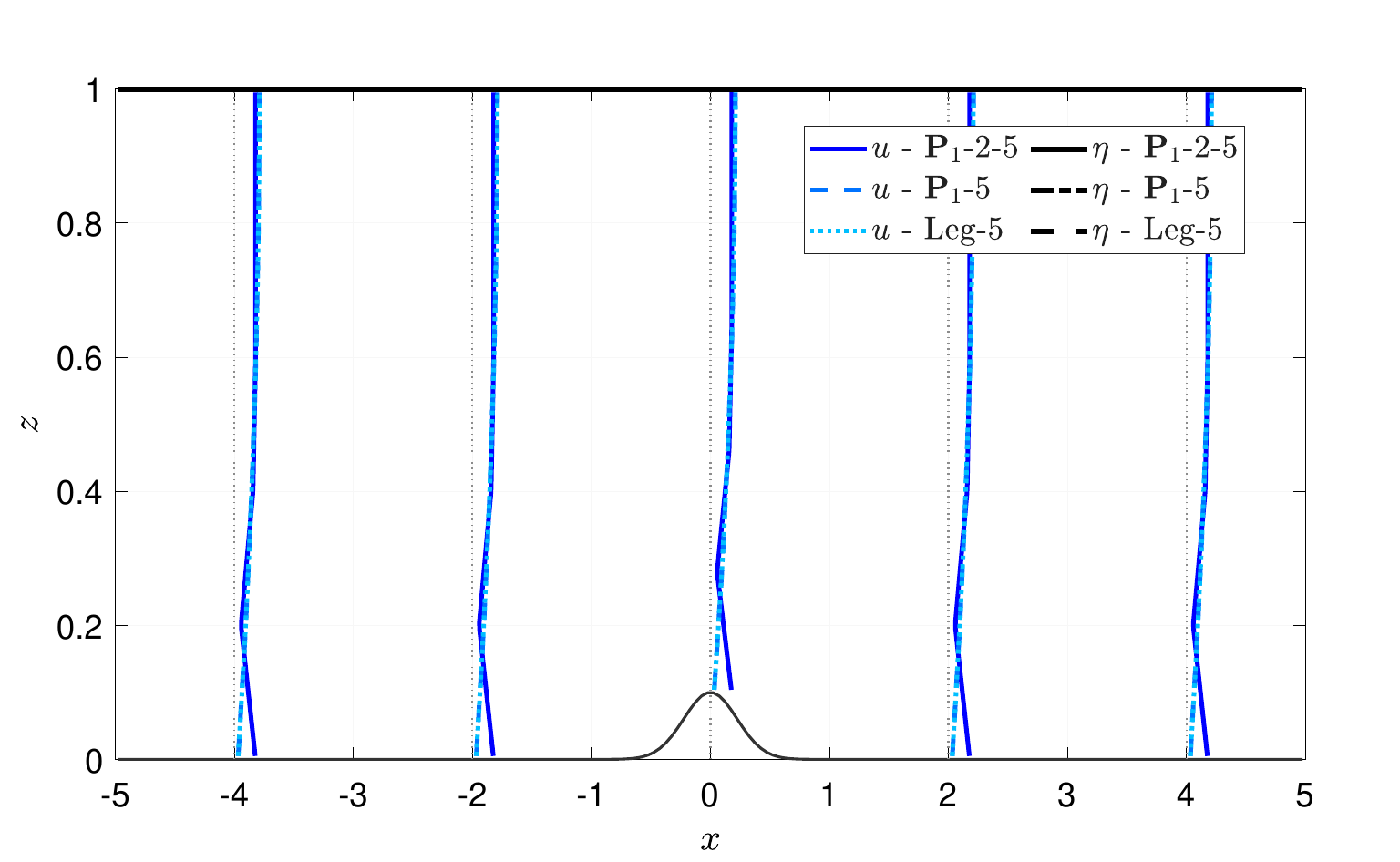}
	\includegraphics[width=0.47\textwidth]{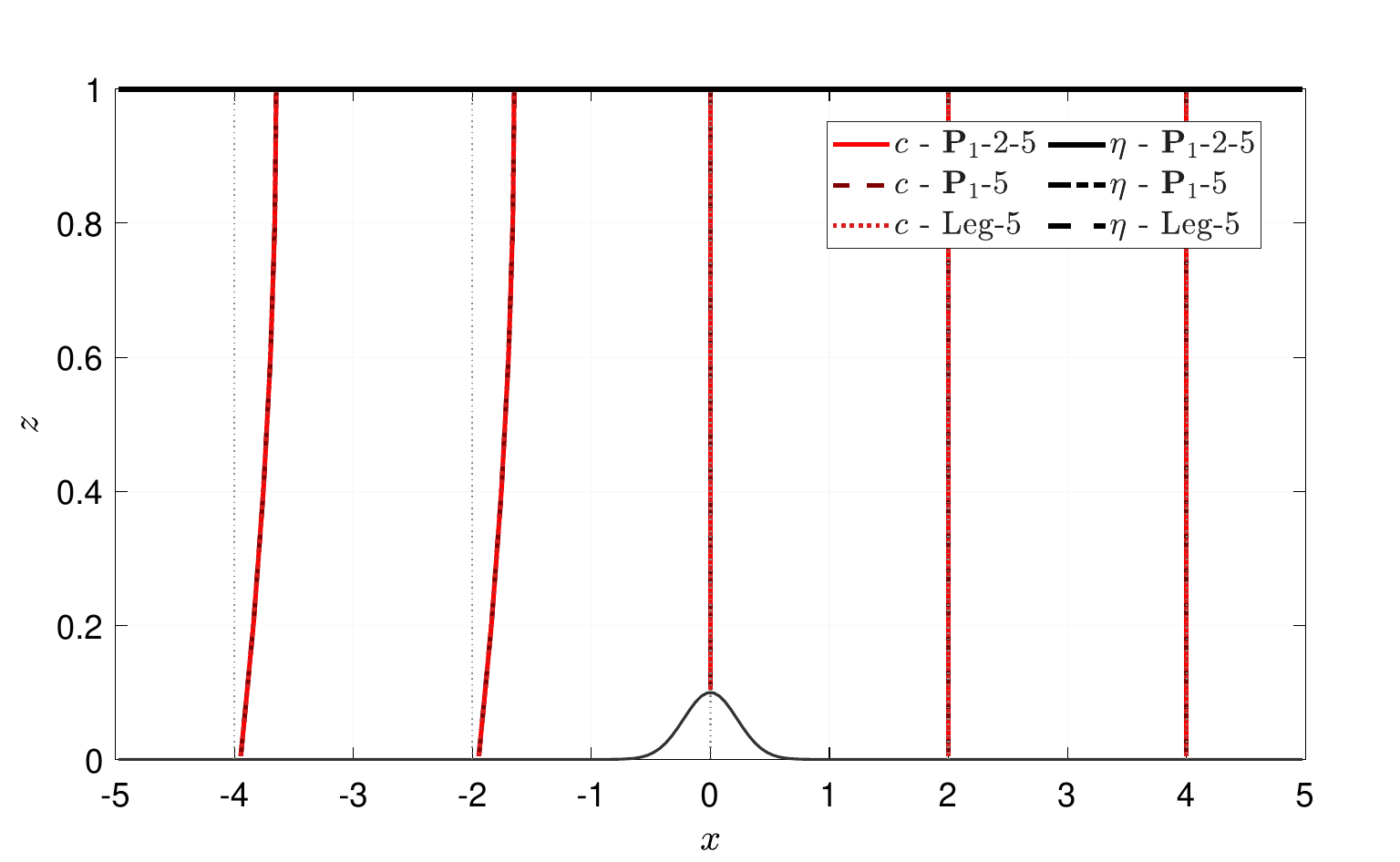}
	\caption{Initial conditions for the velocity and sediment concentration fields for all the considered function bases.}
	\label{fig:test4_ini}
\end{figure}

The initial condition of the moment coefficients is computed by the usual projection onto the corresponding basis (Leg-5, $\mathbb{P}_1$-5 or $\mathbb{P}_1$-2,5) as
$$
u_j(x,0) = \displaystyle\int_0^1 u_{\mathrm{ref}}(\xi)\,\phi_j(\xi) \, d\xi,\qquad\mbox{and}\qquad c_k(x,0) = \displaystyle\int_0^1 c_{\mathrm{ref}}(x,\xi)\,\phi_k (\xi)\, d\xi,\qquad 
$$
for $j=1,\dots,N$ and $k=1,\dots,M$. 

Initial conditions for the velocity and sediment concentration fields are shown in Figure \ref{fig:test4_ini}.  Notice that we compare three configurations for the velocity field: a $\mathbb{P}_1$ reconstruction restricted to the subspace spanned by the first three modes ($N=2$); a $\mathbb{P}_1$ reconstruction with the full set of modes ($N=5$); and a Legendre reconstruction (\(N=5\)). Again, five modes are considered for the sediment concentration approximation in all cases.  We see then that the three reconstructions of the initial concentration profile are identical, while we see small differences, namely close to the bottom, between the bases Leg-5 and $\mathbb{P}_1$-5, and the basis $\mathbb{P}_1$-2-5. This is expected since the approximation space for the velocity field has been restricted.

\begin{figure}[!ht]
	\centering
	\includegraphics[width=0.47\textwidth]{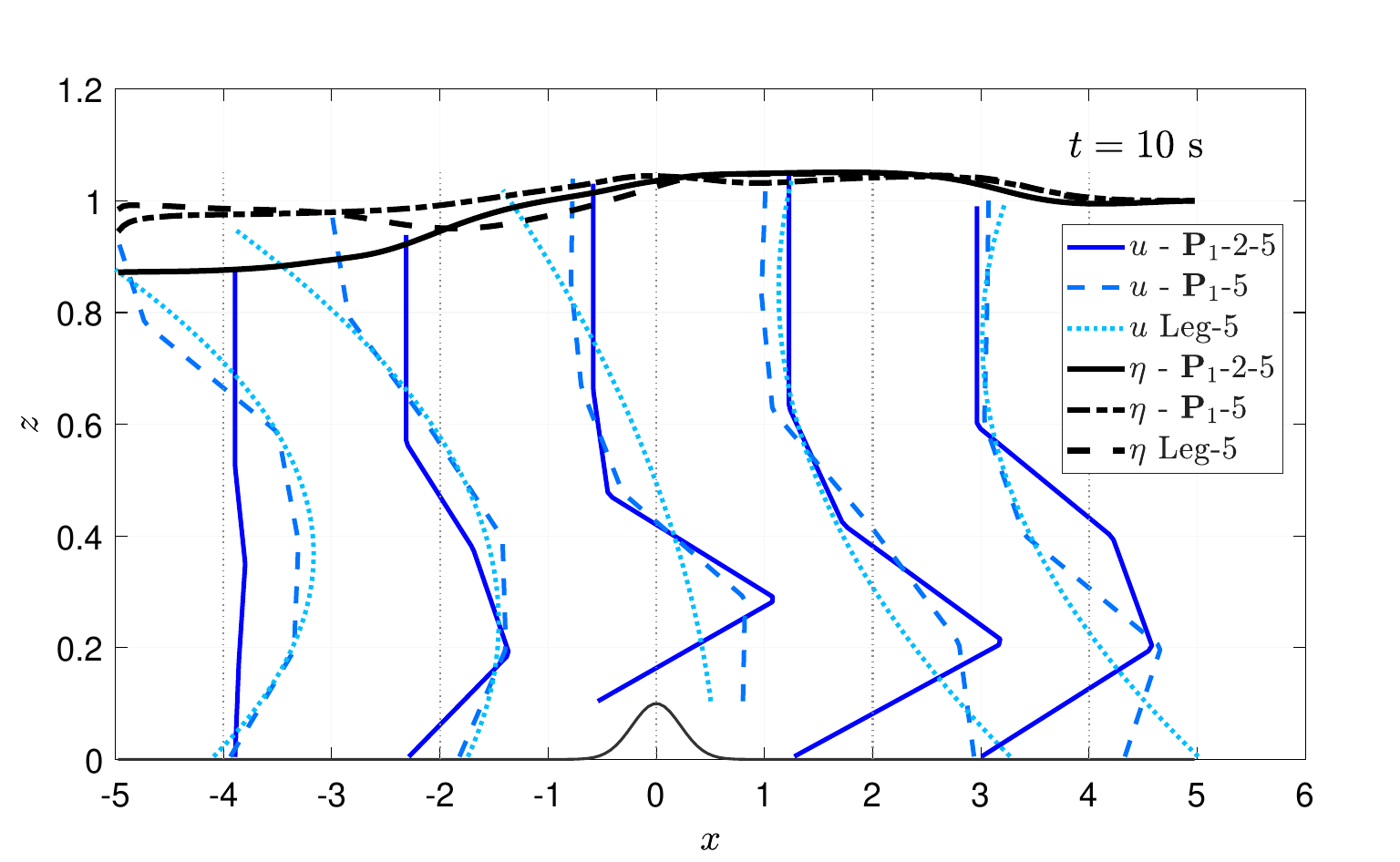}
	\includegraphics[width=0.47\textwidth]{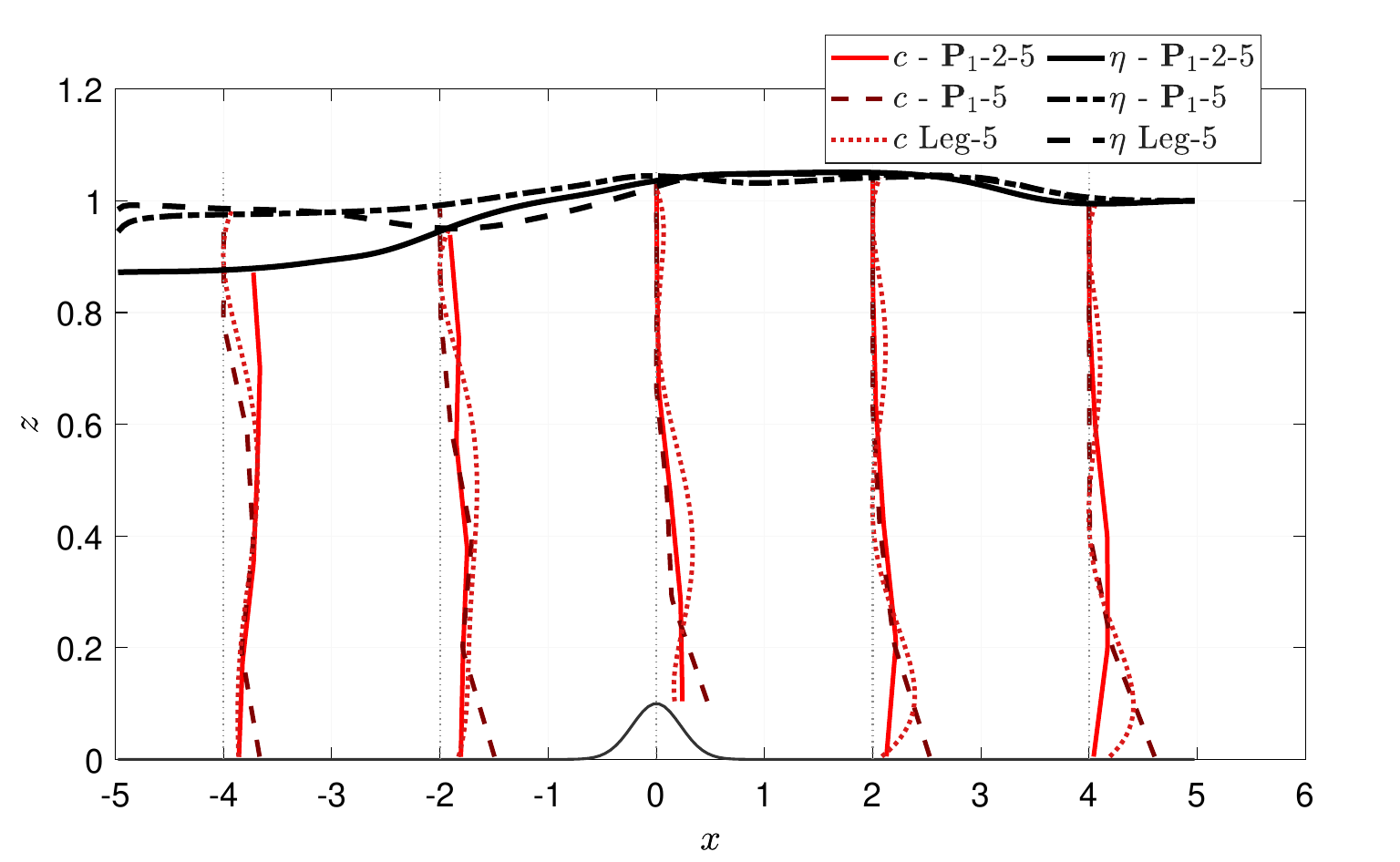}\\
	\includegraphics[width=0.47\textwidth]{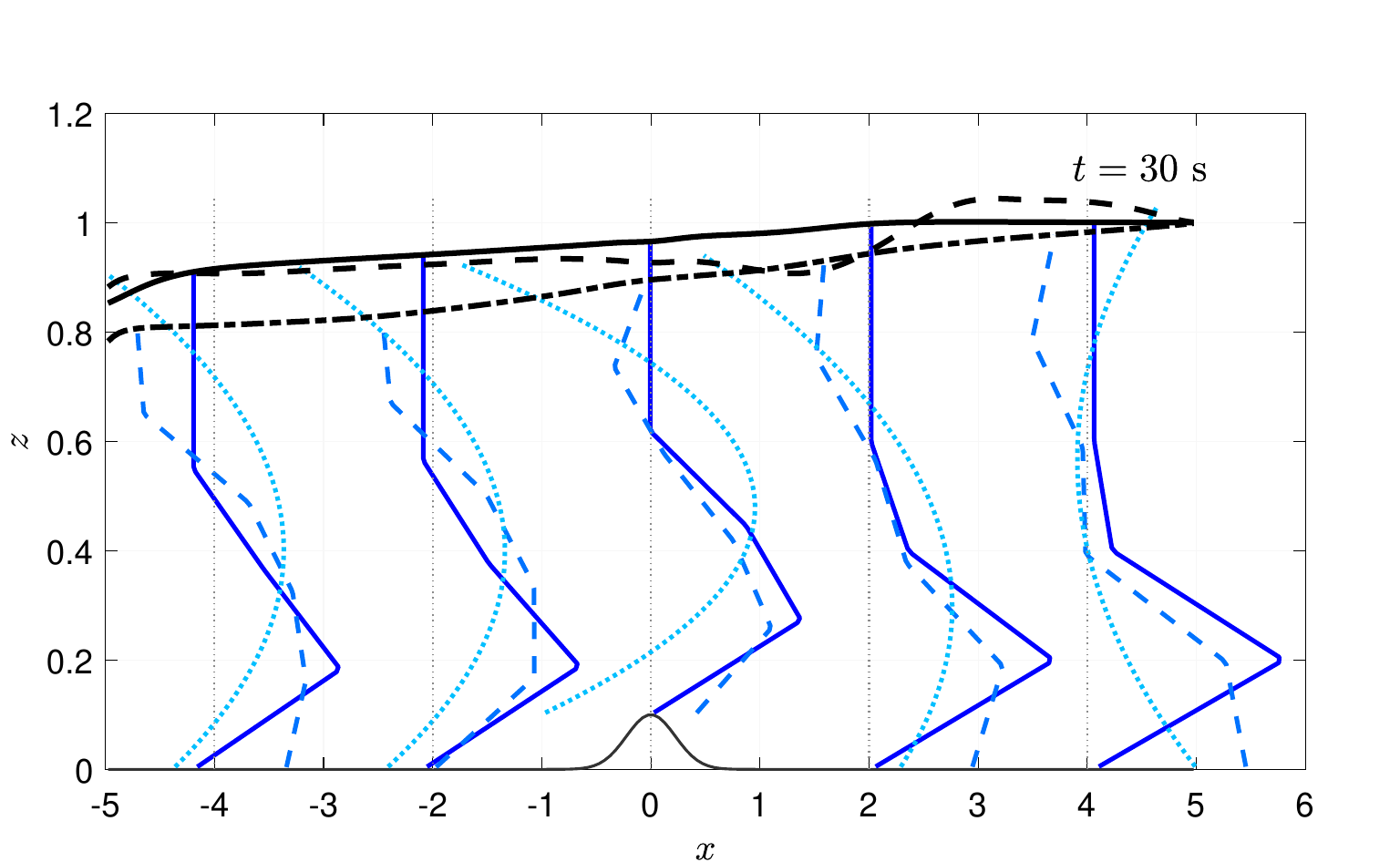}
	\includegraphics[width=0.47\textwidth]{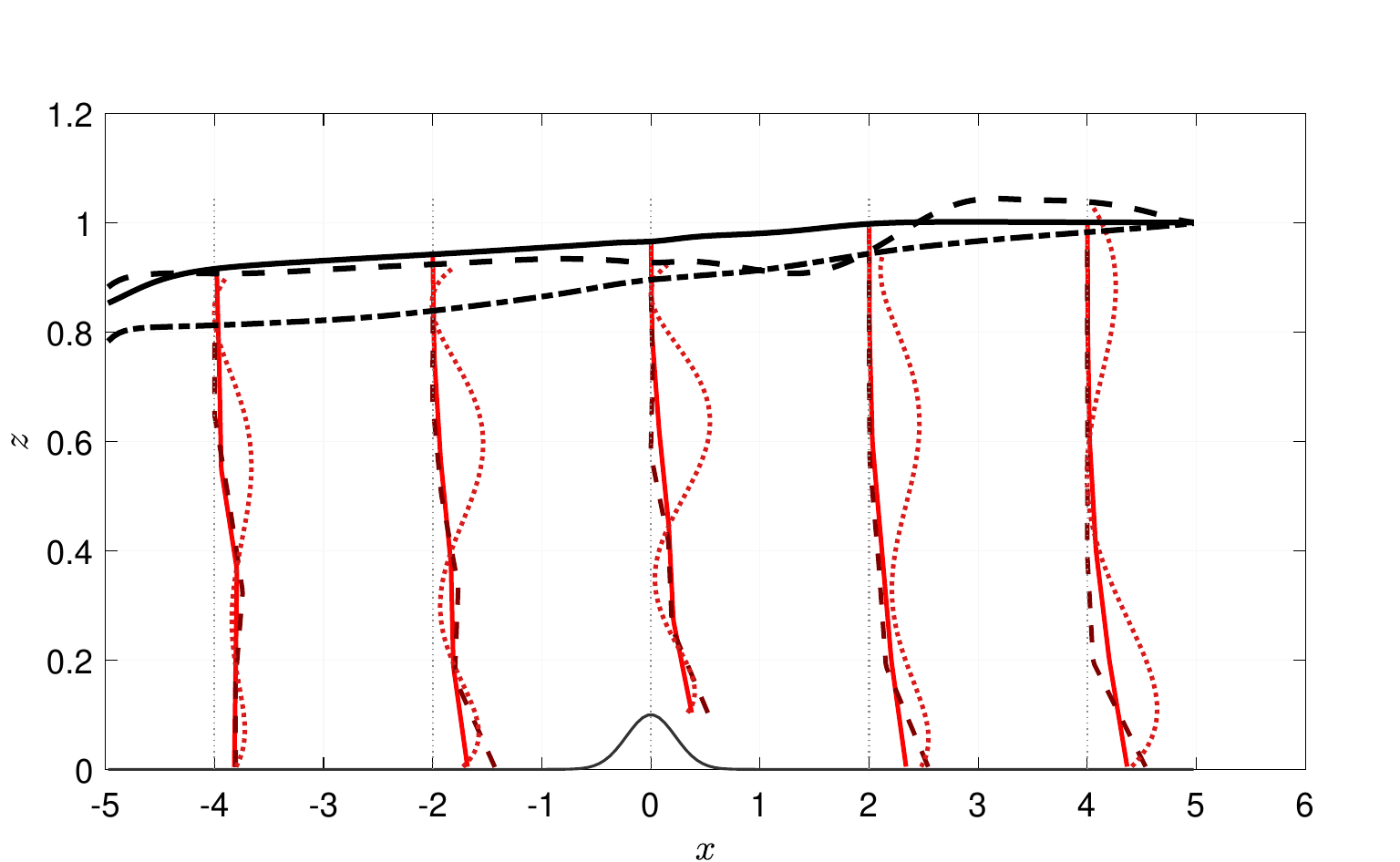}\\
	\includegraphics[width=0.47\textwidth]{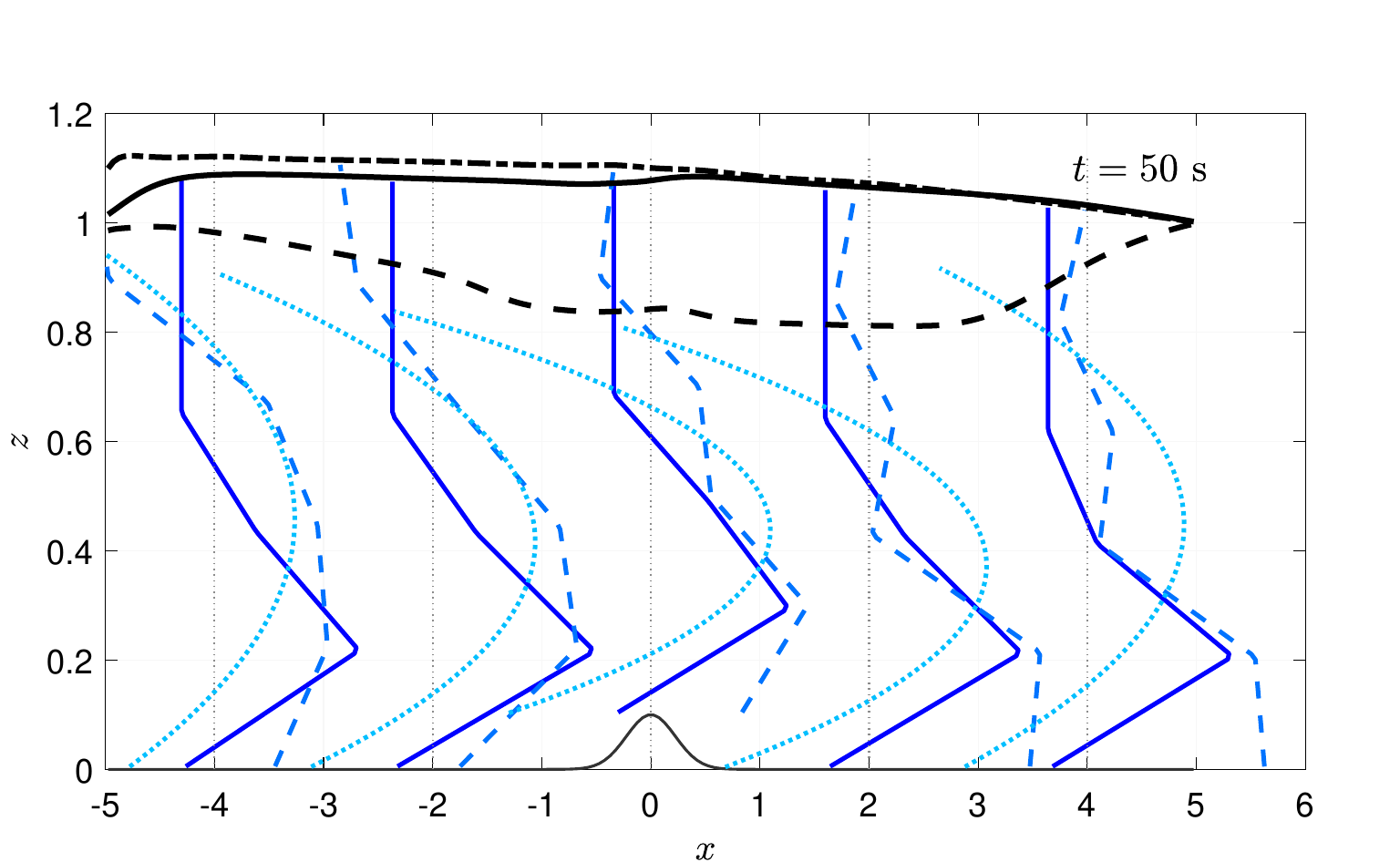}
	\includegraphics[width=0.47\textwidth]{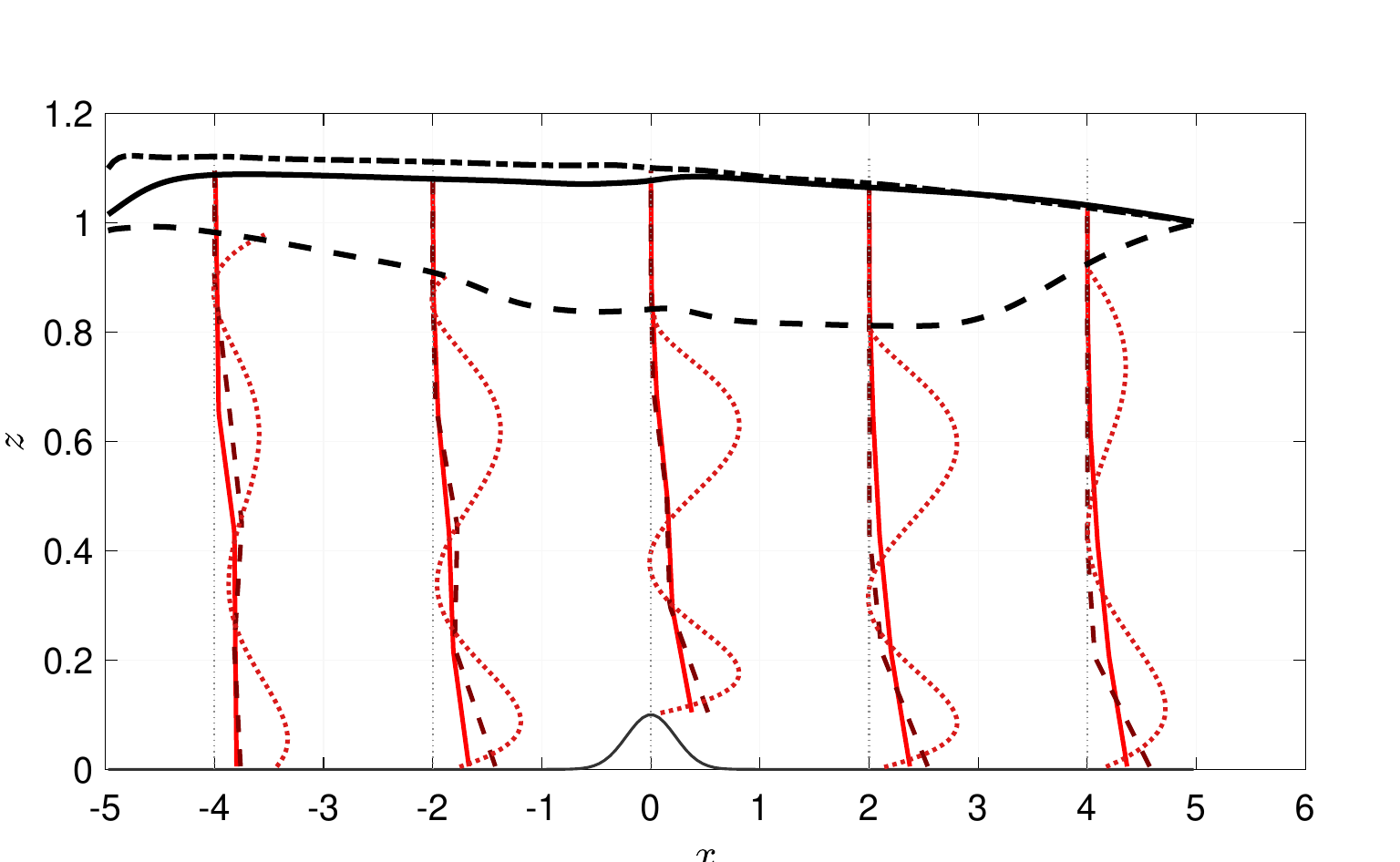}\\
	\includegraphics[width=0.47\textwidth]{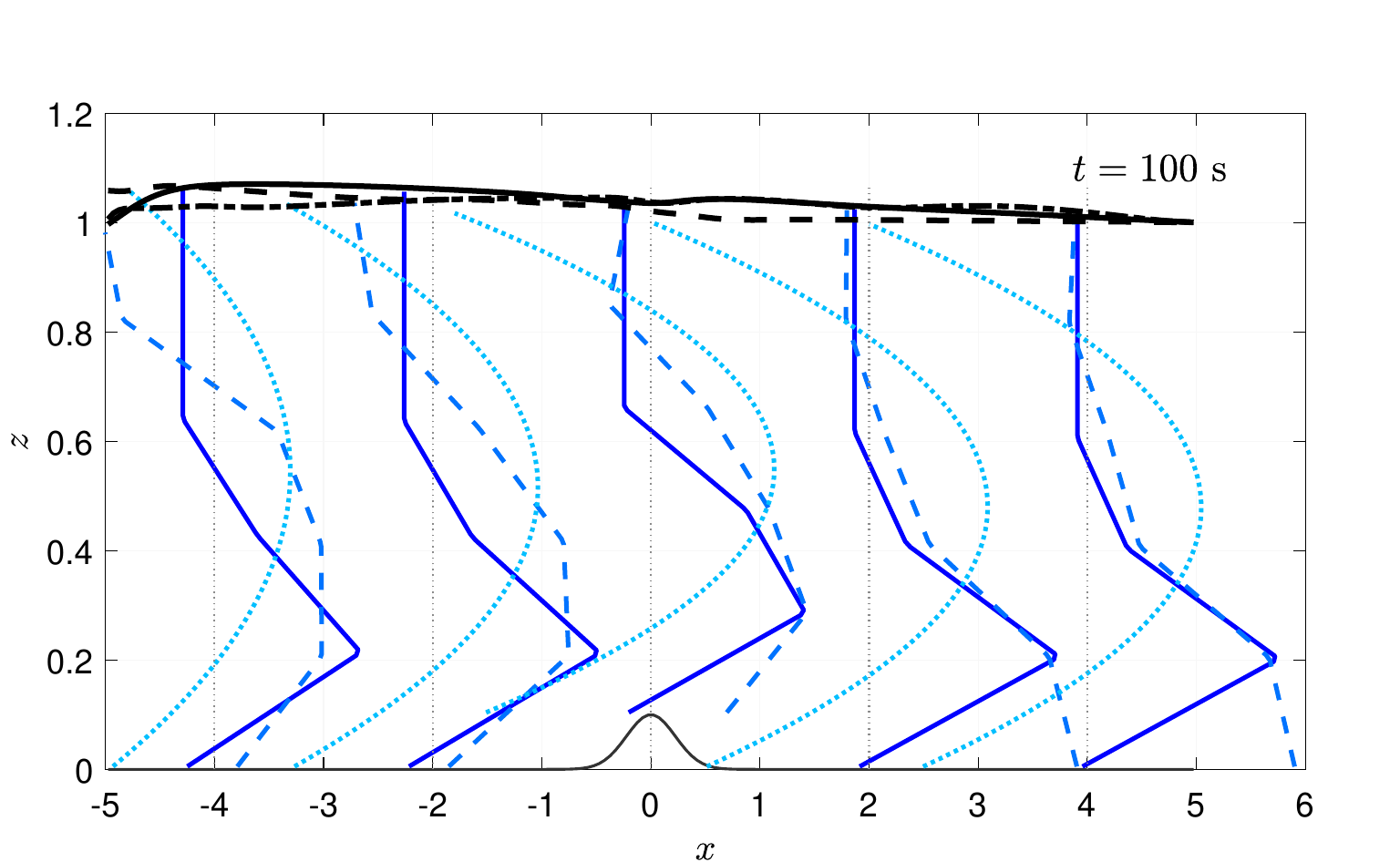}
	\includegraphics[width=0.47\textwidth]{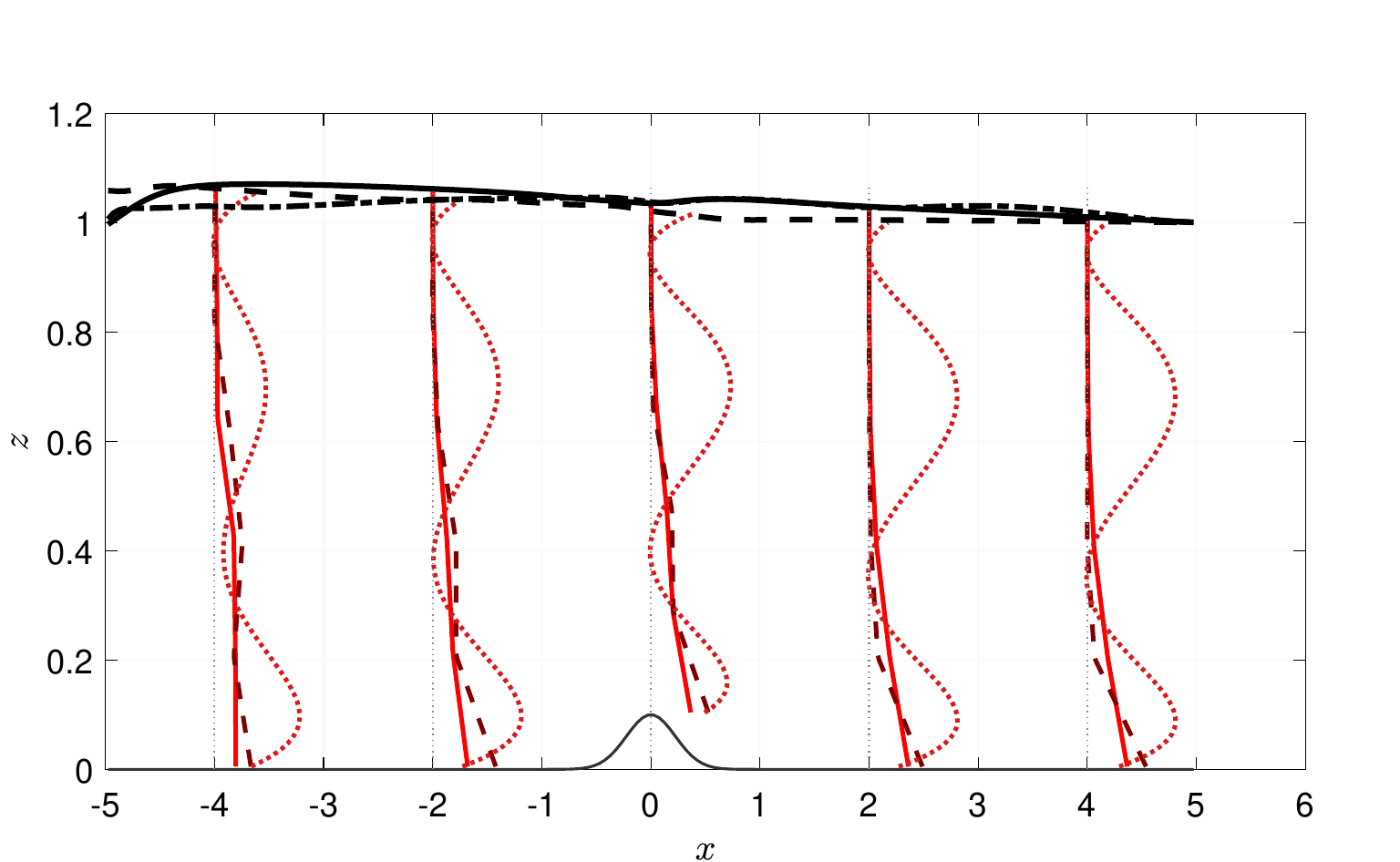}
	\caption{Time evolution of the reconstructed velocity (left) and concentration (right) profiles over the
		smooth bump for times $t=10,30,50,100$ (top to bottom). The $\mathbb{P}_1$ solutions with \(N=2\) and \(N=5\) are compared
		with the Legendre solution with \(N=5\). The concentration is reconstructed with \(M=5\) in all
		simulations. Hence, the differences between the curves are mainly induced by
		the basis choice and by the number of velocity moments \(N\), which affect
		the coupled transport of sediment}
	\label{fig:bump_velocity_concentration_evolution}
\end{figure}

Figure~\ref{fig:bump_velocity_concentration_evolution} shows the evolution of the velocity (left-hand side) and sediment concentration (right-hand side)
profiles. Looking at the velocity field, as the flow interacts with the bottom
bump, a non-trivial vertical shear develops. The $\mathbb{P}_1$ approximation with
\(N=2\) captures the dominant low-order structure but, as expected, the profiles are more
restricted because only three velocity modes are retained. In particular, the moment coefficients $u_k=0$ for $k=3,4,5$, which makes that the corresponding functions basis do not influence the vertical velocity profiles. For this particular choice of the function basis $\mathbb{P}_1$, it makes appear a straight line in the profiles for $\xi>0.6$. Increasing the number of $\mathbb{P}_1$ velocity moments to \(N=5\) gives a richer vertical structure and a smoother evolution of the profile shape. The Legendre approximation with \(N=5\) also captures the vertical shear.

We focus now on the sediment concentration profiles. The sediment initially localized on the left is transported through the computational domain while the settling term
modifies its vertical distribution. Since \(M=5\) is fixed in all simulations,
the differences between the profiles are mainly caused by the different
basis. The $\mathbb{P}_1$-based solutions remain comparatively
localized and regular. The Legendre solution, due to the global character of the
basis, exhibits more pronounced vertical oscillations in the reconstructed
concentration profiles, especially at intermediate and final times. Thus, $\mathbb{P}_1$-based models produce a qualitative physically correct shape of the vertical profile of the sediment concentration, in contrast to polynomial bases that lead to physically meaningless results due to the inherent oscillatory nature of polynomials.

To complement the previous qualitative comparison, we also measured the CPU cost
of the vertical sedimentation update. This step corresponds to the local
vertical evolution associated with the settling term, which is applied after the
hyperbolic update. The measurements were performed for uniform $\mathbb{P}_1$ bases in
two configurations: $\mathbb{P}_1$-2 and $\mathbb{P}_1$-5.

The total CPU time increases from approximately \(2.70\) s for
$\mathbb{P}_1$-2 to approximately \(7.05\) s for
$\mathbb{P}_1$-5, corresponding to a factor of about \(2.61\). However,
the time spent in the sedimentation update increases only from approximately
\(1.88\) s to \(2.84\) s, that is, by a factor of about \(1.51\). Consequently,
the relative contribution of the sedimentation step decreases from about
\(69.7\%\) of the total CPU time for $\mathbb{P}_1$-2 to about
\(40.3\%\) for $\mathbb{P}_1$-5. These measurements indicate that, when the
number of moments is high, the cost is mainly associated with
the hyperbolic coupled moment system, including the evaluation of the moment
fluxes and non-conservative terms. In addition, the cost associated to this step could be reduced using a coarser grid since $N_z=100$ has been considered in all cases.  

\section{Conclusions}\label{se:conclusions}
In this work, a general framework has been derived for hydrostatic shallow moment models for suspended sediment transport with a vertically varying concentration profile. The variables (horizontal velocity and concentration) are expanded in an arbitrary functional basis that is adapted to the specific application, which is subsequently orthogonalised numerically via the Gram-Schmidt algorithm. This approach yields a comprehensive moment model that accounts for the vertical structure of both the velocity and the suspended sediment concentration, avoiding the complex systems that are obtained when using non-orthogonal bases. Concretely, we consider a continuous piecewise linear $\mathbb{P}_1$ finite-element basis. Moreover, non-uniform vertical partitions are allowed in order to obtain an enhanced resolution of the profiles in specific regions of interest. 

The numerical results demonstrate the accuracy and efficiency of this approach. For velocity profiles with boundary layers, the proposed $\mathbb{P}_1$ basis provides accurate approximations with fewer modes. The improvement is even more significant for suspended sediment concentration in sedimentation processes, whose profiles frequently exhibit very sharp vertical gradients. In such cases, using global Legendre polynomial bases produces vertical profiles with severe oscillations that lead to physically meaningless results. On the contrary, the piecewise linear reconstructions give as result qualitatively accurate profiles that are in agreement with the expected concentration distributions observed in sedimentation processes.

Furthermore, the proposed formulation allows for an asymmetric framework where each variable can be projected onto a different subspace. This feature is particularly useful here since the sediment concentration often requires a richer approximation space than the fluid velocity. Therefore, the proposed framework provides a highly flexible and computationally efficient way to adapt the vertical discretization to the specific physical behaviour of each variable in coupled geophysical flows. In the future, the application of these techniques to other physical contexts where standard moment models fail due to the complex shape of the vertical profiles will be studied.

\section*{Acknowledgements}
This research has been partially supported by the European Union - NextGenerationEU program and by grant PID2022-137637NB-C2(1/2) and PID2025-170651NB-C2(1/2) funded by  MCIN/AEI/10.13039/501100011033 and ``ERDF A way of making Europe''. J. Garres-Díaz thanks IMUS-Maria de Maeztu grant CEX2024-001517-M - Apoyo a Unidades de Excelencia María de Maeztu for supporting this research, funded by MICIU/AEI/10.13039/501100011033 and the applied research and innovation project "Development of mathematical tools for transfer" (Mathware) at IMUS, co-financed by the EU - Ministry of Finance and Public Administration - European Funds - Junta de Andalucía -Consejería de Universidad, Investigación e Innovación. I. Gómez-Bueno and T. Morales de Luna thank project PCI2024-155061-2 funded by MCIN/AEI/10.13039/501100011033 and ``ERDF A way of making Europe''.

\bibliographystyle{plain}

\begin{thebibliography}{10}

\bibitem{Bauerle2025}
M.~Bauerle, A-J. Christlieb, M.~Ding, and J.~Huang.
\newblock On the rotational invariance and hyperbolicity of shallow water
  moment equations in two dimensions.
\newblock {\em SIAM Journal on Mathematical Analysis}, 57(1):1039--1085, 2025.

\bibitem{bradford1999hydrodynamics}
S.~F. Bradford and N.~D. Katopodes.
\newblock Hydrodynamics of turbid underflows. {I}: Formulation and numerical
  analysis.
\newblock {\em Journal of Hydraulic Engineering}, 125(10):1006--1015, 1999.

\bibitem{Buerger2019}
R.~B\"urger, E.D. Fern\'andez-Nieto, and V.~Osores.
\newblock A dynamic multilayer shallow water model for polydisperse
  sedimentation.
\newblock {\em ESAIM: M2AN}, 53(4):1391--1432, 2019.

\bibitem{Burger2007}
R.~B\"urger, A.~Garc\'ia, K.~H. Karlsen, and J.~D. Towers.
\newblock A family of numerical schemes for kinematic flows with discontinuous
  flux.
\newblock {\em Journal of Engineering Mathematics}, 60(3-4):387--425, 2007.

\bibitem{CanteroChinchilla2019}
F.~N. Cantero-Chinchilla, O.~Castro-Orgaz, and A.~A. Khan.
\newblock Vertically-averaged and moment equations for flow and sediment
  transport.
\newblock {\em Advances in Water Resources}, 132:103387, oct 2019.

\bibitem{Castro2017}
M.~J. Castro, T.~Morales~de Luna, and C.~Par{\'e}s.
\newblock Well-{Balanced} {Schemes} and {Path}-{Conservative} {Numerical}
  {Methods}.
\newblock In R\'emi Abgrall {and} Chi-Wang Shu, editor, {\em Handbook of
  {Numerical} {Analysis}}, volume~18 of {\em Handbook of {Numerical} {Methods}
  for {Hyperbolic} {Problems Applied} and {Modern} {Issues}}, pages 131--175.
  Elsevier, 2017.
\newblock DOI: 10.1016/bs.hna.2016.10.002.

\bibitem{CastroDiaz2012}
M.~J. Castro~D{\'\i}az and E.~D. Fern{\'a}ndez-Nieto.
\newblock A class of computationally fast first order finite volume solvers:
  {PVM} methods.
\newblock {\em {SIAM} Journal on Scientific Computing}, 34(4):A2173--A2196, jan
  2012.

\bibitem{SaintVenant1871}
A.~J.~C.~.B. de~Saint-Venant.
\newblock Th{\'e}orie du mouvement non permanent des eaux, avec application aux
  crues des rivi\`{e}res et \`{a} l'introduction des mar{\'e}es dans leur lit.
\newblock {\em C. R. Acad{\'e}mie des Sciences, Paris}, 73:147--154, 1871.

\bibitem{Dey2014}
S.~Dey.
\newblock {\em Fluvial Hydrodynamics}.
\newblock Springer-Verlag GmbH, 2014.

\bibitem{FernandezNieto2013}
E.~D. Fern{\'a}ndez-Nieto, E.~H. Kon{\'e}, T.~Morales~de Luna, and
  R.~B{\"u}rger.
\newblock A multilayer shallow water system for polydisperse sedimentation.
\newblock {\em Journal of Computational Physics}, 238:281--314, April 2013.

\bibitem{FernandezNieto2013a}
E.~D. Fern{\'{a}}ndez-Nieto, E.~H. Kon{\'{e}}, and T.~Chac{\'{o}}n Rebollo.
\newblock A multilayer method for the hydrostatic {N}avier-{S}tokes equations:
  A particular weak solution.
\newblock {\em Journal of Scientific Computing}, 60(2):408--437, nov 2013.

\bibitem{GarresDiaz2021b}
J.~Garres-D{\'{i}}az, M.~J. Castro~D{\'{i}}az, J.~Koellermeier, and
  T.~Morales~de Luna.
\newblock Shallow water moment models for bedload transport problems.
\newblock {\em Communications in Computational Physics}, 30(3):903--941, jun
  2021.

\bibitem{GarresDiaz2023}
J.~Garres-D{\'{i}}az, C.~Escalante, T.~{Morales de Luna}, and M.J. {Castro
  D{\'{i}}az}.
\newblock A general vertical decomposition of {E}uler equations:
  {M}ultilayer-moment models.
\newblock {\em Applied Numerical Mathematics}, 183:236--262, jan 2023.

\bibitem{Gerbeau2001}
J.-F. Gerbeau and B.~Perthame.
\newblock Derivation of viscous {S}aint-{V}enant system for laminar shallow
  water; numerical validation.
\newblock {\em Discrete Contin. Dyn. Syst. Ser. B}, 1(1):89--102, 2001.

\bibitem{GonzalezAguirre2020}
J.~C. Gonz\'alez-Aguirre, M.~J. Castro, and T.~Morales de~Luna.
\newblock A robust model for rapidly varying flows over movable bottom with
  suspended and bedload transport: {M}odelling and numerical approach.
\newblock {\em Advances in {W}ater {R}esources}, 140:103575, 2020.

\bibitem{Guo1999}
Q.-C. Guo and Y.-C. Jin.
\newblock Modeling sediment transport using depth-averaged and moment
  equations.
\newblock {\em Journal of Hydraulic Engineering}, 125(12):1262--1269, December
  1999.

\bibitem{khan2005numerical}
S.~M. Khan, J.~Imran, S.~Bradford, and J.~Syvitski.
\newblock Numerical modeling of hyperpycnal plume.
\newblock {\em Marine Geology}, 222-223:193--211, November 2005.

\bibitem{Koellermeier2022}
J.~Koellermeier and E.~Pimentel-Garc\'ia.
\newblock Steady states and well-balanced schemes for shallow water moment
  equations with topography.
\newblock {\em Applied Mathematics and Computation}, 427:127166, August 2022.

\bibitem{Koellermeier2020}
J.~Koellermeier and E.~Pimentel-Garc\'ia.
\newblock Steady states and well-balanced schemes for shallow water moment
  equations with topography.
\newblock {\em Applied Mathematics and Computation}, 427:127166, 2022.

\bibitem{koellermeier2020analysis}
J.~Koellermeier and M.~Rominger.
\newblock Analysis and numerical simulation of hyperbolic shallow water moment
  equations.
\newblock {\em Communications in Computational Physics},
  28(3):1038--1084, 2020.

\bibitem{Kowalski2018}
J.~Kowalski and M.~Torrilhon.
\newblock Moment approximations and model cascades for shallow flow.
\newblock {\em Communications in Computational Physics}, 25(3):669--702, 2018.

\bibitem{MoralesdeLuna2017}
T.~Morales~de Luna, E.~D. Fern\'{a}ndez~Nieto, and M.~J. Castro~D\'{i}az.
\newblock Derivation of a multilayer approach to model suspended sediment
  transport: Application to hyperpycnal and hypopycnal plumes.
\newblock {\em Communications in Computational Physics}, 22(5):1439--1485,
  2017.

\bibitem{Pares2006}
C.~Par{\'e}s.
\newblock Numerical methods for nonconservative hyperbolic systems: a
  theoretical framework.
\newblock {\em {SIAM} Journal on Numerical Analysis}, 44(1):300--321, jan 2006.

\bibitem{parker1986self}
G.~Parker, Y.~Fukushima, and H.~M. Pantin.
\newblock Self-accelerating turbidity currents.
\newblock {\em Journal of Fluid Mechanics}, 171:145--181, 1986.

\bibitem{Burger2025}
{R. Bürger}, {E.~D. Fern{\'a}ndez-Nieto}, {J. Garres-D{\'i}az}, and {J. Moya}.
\newblock Well-balanced physics-based finite volume schemes for
  {S}aint-{V}enant-{E}xner-type models of sediment transport.
\newblock {\em Advances in Water Resources}, 206:105178, December 2025.

\bibitem{Koellermeier2025}
U.~Scholz and J.~Koellermeier.
\newblock Spline shallow water moment equations.
\newblock arXiv preprint arXiv:2508.02714, 2025.

\bibitem{Steffler1993}
P.~M. Steffler and J.~Y.-Chung.
\newblock Depth averaged and moment equations for moderately shallow free
  surface flow.
\newblock {\em Journal of Hydraulic Research}, 31(1):5--17, jan 1993.

\end{thebibliography}

\end{document}